\documentclass[review,3p,square,11pt]{elsarticle}

\usepackage{amssymb}
\usepackage{amsmath}
\usepackage{amsthm}
\usepackage{stmaryrd}
\usepackage{bbm}
\usepackage{enumitem}

\usepackage{fancyhdr}
\usepackage{hyperref}
\hypersetup{%
  pdftitle={BSDEs},%
  pdfsubject={BSDEs},%
  pdfauthor={ShengJun FAN, Ying Hu},%
  pdfkeywords={BSDEs},%
  pdfstartview=FitH,%
  CJKbookmarks=true,%
  bookmarksnumbered=true,%
  bookmarksopen=true,%
  colorlinks=true, linkcolor=blue, urlcolor=blue, citecolor=blue, %
}

\usepackage{cleveref}
\crefname{thm}{Theorem}{Theorems}
\crefname{pro}{Proposition}{Propositions}
\crefname{lem}{Lemma}{Lemmas}
\crefname{rmk}{Remark}{Remarks}
\crefname{cor}{Corollary}{Corollaries}
\crefname{dfn}{Definition}{Definitions}
\crefname{ex}{Example}{Examples}
\crefname{section}{Section}{Sections}
\crefname{subsection}{Subsection}{Subsections}

\newcommand{\eps}{\varepsilon}

\newcommand{\To}{\rightarrow}
\newcommand{\as}{{\rm d}\mathbb{P}\times{\rm d} t-a.e.}
\newcommand{\ass}{{\rm d}\mathbb{P}\times{\rm d} s-a.e.}
\newcommand{\ps}{\mathbb{P}-a.s.}

\newcommand{\F}{\mathcal{F}}
\newcommand{\E}{\mathbb{E}}

\newcommand{\bs}{\mathbf{S}}

\newcommand{\T}{[0,T]}

\newcommand{\R}{{\mathbb R}}

\newcommand{\RE}{\forall}

\newtheorem{thm}{Theorem}[section]

\newtheorem{pro}[thm]{Proposition}
\newtheorem{rmk}[thm]{Remark}
\newtheorem{cor}[thm]{Corollary}

\newtheorem{ex}[thm]{Example}

\journal{arXiv}

\begin{document}
\begin{frontmatter}

\title{{Solvability of one-dimensional sub-quadratic BSDEs with unbounded stochastic coefficients}\tnoteref{found}}
\tnotetext[found]{This work is supported by the National Natural Science Foundation of China (Nos. 12171471 and 12631016).\vspace{0.2cm}}


\author[Fan,Ji]{Haibo Ji} \ead{jhb\_2010@163.com}
\author[Fan]{Shengjun Fan\corref{cor}} \ead{shengjunfan@cumt.edu.cn}
\vspace{-0.5cm}

\affiliation[Fan]{organization={School of Mathematics, China University of Mining and Technology},
            city={Xuzhou 221116},
            country={China}}

\affiliation[Ji]{organization={School of Mathematics and Physics, Suqian University},
            city={Suqian  223800},
            country={China}}\vspace{-0.8cm}

\cortext[cor]{Corresponding author\vspace{0.2cm}}

\begin{abstract}
This study is devoted to solvability of the unbounded solutions of a class of one-dimensional sub-quadratic backward stochastic differential equations (BSDEs) with unbounded stochastic coefficients and general time horizons. In a general sub-exponential moment integrability space linked with the stochastic coefficients, we first establish existence of an unbounded solution of the preceding BSDE when the generator $g$ has a one-sided linear growth in the state variable $y$ and a sub-quadratic growth in the state variable $z$, and then prove a comparison theorem of the unbounded solutions of such two BSDEs when either of generators satisfies an extended convexity/concavity condition in $(y,z)$ with unbounded stochastic coefficients, yielding naturally the uniqueness of the solution. Particularly, the unbounded stochastic coefficient before the state variable $y$ is supposed to take values in $\R$ such that the required integrability condition of the terminal value may be weaker than that under the usual situations. These findings incorporate and generalize some corresponding known results. \vspace{0.2cm}
\end{abstract}

\begin{keyword}
Backward stochastic differential equation \sep Existence and uniqueness \sep \\
\hspace*{2.1cm}Sub-quadratic growth\sep Comparison theorem \sep Unbounded stochastic coefficients \vspace{0.2cm}

\MSC[2010] 60H10\vspace{0.2cm}
\end{keyword}

\end{frontmatter}
\vspace{-0.5cm}

\section{Introduction}
\label{sec:1-Introduction}
\setcounter{equation}{0}

We fix an extended positive real $T$ and a positive integer $d$, and let $(\Omega, \mathcal{F}, \mathbb{P})$ be a complete probability space, $(B_t)_{t\in[0,T ]}$ be an $\mathbb{R}^d$-valued standard Brownian motion defined on this space, and $(\mathcal{F}_t)_{t\in[0,T ]}$
the natural filtration generated by $(B_t)_{t\geq0}$ and augmented by all P-null sets
of $\mathcal{F}$. Any progressive measurability and adaptability of processes will be associated to this
filtration, and all equalities and inequalities between random elements will be understood to hold $\ps$. We consider the following backward stochastic differential equation (BSDE for short):
\begin{equation}\label{eq:1.1}
Y_t=\xi+\int_{t}^{T} g(s,Y_s,Z_s)\mathrm{d}s-\int_{t}^{T}Z_s\cdot \mathrm{d}B_s,\ \ t\in\T,
\end{equation}
where the terminal condition $\xi$ is an $\mathcal{F}_T $-measurable real-valued random variable, and the
generator $g(\omega,t,y,z):\Omega\times[0,T]\times\mathbb{R}\times\mathbb{R}^d\rightarrow\mathbb{R}$ is $(\mathcal{F}_t)$-progressively measurable for any $(y,z)$ and continuous in $(y,z)$. The triple $(\xi,T,g)$ is called the parameters of BSDE \eqref{eq:1.1}, denoted by BSDE $(\xi,T,g)$. The pair $(Y_t,Z_t)_{t\in[0,T]}$ of $(\mathcal{F}_t)$-progressively measurable processes taking values in $\R\times\R^d$ and satisfying \eqref{eq:1.1} is called an adapted solution of BSDE \eqref{eq:1.1}.\vspace{0.2cm}

It is well known that \citet{Bi} first proposed the linear BSDEs, which were regarded as the adjoint processes in the maximum principle of stochastic control, and that nonlinear BSDEs were initially studied by \citet{PP} in 1990, which verified an existence and uniqueness theorem for the square-integrable solutions of multidimensional BSDEs with uniformly Lipschitz continuous generators, and was called
the ``Founder Paper" of the BSDE theory. Since then, BSDEs have received widespread attention and gained more applications. They have increasingly been proved to be a powerful tool in various fields such as nonlinear mathematical expectation, partial differential equation, stochastic control, mathematical finance and so on. More details can be seen in \cite{P1,P2,EPQ,DP,Pa,K,HIM,T,FWY,FanHuTang2026MF} among others.\vspace{0.2cm}

In this paper, we study the unbounded adapted solution of BSDE \eqref{eq:1.1} with a sub-quadratic generator $g$, which satisfies, roughly speaking, that $\as$, for each $(y,z)\in \R\times\R^d$,
\begin{equation}\label{eq:1.2}
\mathrm{sgn}(y)g(\omega,t,y,z)\leq f_t(\omega)+\beta_{t}(\omega)|y|+\gamma_t(\omega)|z|^\delta,
\end{equation}
where $(f_t)_{t\in[0,T]}$, $(\beta_t)_{t\in[0,T]}$ and $(\gamma_t)_{t\in[0,T]}$ are three $(\mathcal{F}_t$)-progressively measurable real-valued processes. We would like to emphasize that they are all possibly unbounded and $\beta_\cdot$ can take negative values. If $g$ satisfies \eqref{eq:1.2}, then we say that $g$ has a one-sided linear growth in $y$, a sub-linear growth in $z$ when $\delta\in (0,1)$, a linear growth in $z$ when $\delta=1$, a quadratic growth in $z$ when $\delta=2$, a super-quadratic growth in $z$ when $\delta>2$, and a sub-quadratic growth in $z$ when $\delta\in (1,2)$. In this paper, we would mainly pay attention to the last case. For narrative convenience, when $\delta\in (1,2)$, we let $\delta^*$ be the conjugate of $\delta$, i.e., $1/\delta+1/\delta^*=1$.\vspace{0.2cm}

When both $\beta_\cdot$ and $\gamma_\cdot$ are either nonnegative constants or deterministic nonnegative functions, many scholars have studied solvability of the adapted solutions of corresponding BSDEs for the cases of $\delta\in (0,1)$, $\delta=1$, $\delta\in(1,2)$, $\delta=2$ and $\delta>2$, respectively. Generally speaking, to obtain existence/uniqueness of the adapted solution, the terminal conditions $\xi$ and $\int_0^T f_t{\rm d}t$ are required to have at least an $L^1$-integrability when $\delta\in (0,1)$, an $L^p$-integrability for some $p>1$ or an $L\exp(\mu\sqrt{\ln(1+L)})$-integrability for some/any $\mu>0$ when $\delta=1$, an $\exp(\mu L^{2\over \delta^*})$-integrability for some/any $\mu>0$ when $\delta\in(1,2)$, and an $\exp(\mu L)$-integrability for some/any $\mu>0$ when $\delta=2$. And, BSDEs with super-quadratic growth generators in $z$ are not solvable in general without the Markovian setting. The reader is referred to \cite{Pa,CW,K,BD,BH,BH1,DHB,BE1,F1,LF,F2,FanHuTang2020CRM,FH1,FanHuTang2026MF} for details. In particular,
by the test function method and the localization technique, \citet{FH1} established existence of the unbounded solution of sub-quadratic BSDE \eqref{eq:1.1} when the terminal time $T<+\infty$, the terminal conditions $\xi$ and $\int_0^T f_t{\rm d}t$ admits only an $\exp(\mu L^{2\over \delta^*})$-integrability for some $\mu>0$, the generator $g$ satisfies condition \eqref{eq:1.2} with $\delta\in (1,2)$ and with $\beta_\cdot$ and $\gamma_\cdot$ being two nonnegative constants $\beta$ and $\gamma$, and it also satisfies the following general growth condition: $\as$, for each $(y,z)\in \R\times\R^d$,
\begin{equation*}
|g(\omega,t,y,z)|\leq f_t(\omega)+\varphi(|y|)+c |z|^2,
\end{equation*}
where $c>0$ is a constant and $\varphi(x):\R_+\To \R_+$ is an nondecreasing continuous function with $\varphi(0)=0$. Meanwhile, by virtue of $\theta$-difference technique, \cite{FH1} also proved uniqueness of the unbounded solution such that the supremum of the value process has a $\cap_{\mu>0}\exp(\mu L^{2\over \delta^*})$-integrability when the terminal conditions admits a same integrability and the generator $g$ further satisfies the following extended convexity/concavity condition related to the sub-quadratic growth: $\as$, for each $(y_i,z_i)\in \R\times\R^d,\, i=1,2$ and each $\theta\in(0,1)$,
\begin{align}\label{eq:1.4}
\mathbf{1}_{\{\Delta_\theta y>0\}}\frac{g(\omega,t,y_1,z_1)-\theta g(\omega,t,y_2,z_2)}{1-\theta}\leq f_t(\omega)+\beta |y_2|+\beta (\Delta_\theta y)^+ +\gamma|\Delta_\theta z|^\delta
\end{align}
or
\begin{align}\label{eq:1.5}
-\mathbf{1}_{\{\Delta_\theta y<0\}}\frac{g(\omega,t,y_1,z_1)-\theta g(\omega,t,y_2,z_2)}{1-\theta}\leq f_t(\omega)+\beta |y_2|+\beta (\Delta_\theta y)^- +\gamma |\Delta_\theta z|^\delta,
\end{align}
where
$$
\Delta_\theta y:=\frac{y_1-\theta y_2}{1-\theta}\ \ {\rm and}\ \ \Delta_\theta z:=\frac{z_1-\theta z_2}{1-\theta}.\vspace{0.2cm}
$$
Furthermore, \citet{GWF} extended the existence and uniqueness results of \cite{FH1} to the case that the terminal time $T$ can be positive infinity and the constant coefficients $\beta$ and $\gamma$ can be deterministic nonnegative functions $\beta(t)$ and $\gamma(t)$ defined on $\T$ and satisfying
$$
\int_0^T \left(\beta(t)+\gamma^{2\over 2-\delta}(t)\right){\rm d}t<+\infty.
$$

On the other hand, when the generator $g$ satisfies \eqref{eq:1.2} with $\delta=1$ and with $\beta_\cdot$ and $\gamma_\cdot$ being two nonnegative stochastic processes such that for some $M>0$
$$
\int_0^T \left(\beta_t(\omega)+\gamma^2_t(\omega)\right){\rm d}t\leq M,
$$
\citet{LiuLiFan2020CAM}, \citet{LiXuFan2021PUQR} and \citet{LiFan2023CSTM} studied solvability of the $L^p~(p>1)$ solutions of general time-interval BSDEs. And, when the generator $g$ satisfies \eqref{eq:1.2} with $\delta=2$, $\beta_\cdot$ being a nonnegative stochastic process such that for some $M>0$
$$
\int_0^T \beta_t(\omega)\,{\rm d}t\leq M,
$$
and $\gamma_\cdot$ being a nonnegative constant, \citet{FL} and \citet{LF} addressed solvability of the bounded solutions of general time-interval BSDEs. To the best of our knowledge, when the generator $g$ satisfies \eqref{eq:1.2} with $\beta_\cdot$ and $\gamma_\cdot$ being general unbounded stochastic processes, most of solvability results on the adapted solutions of BSDEs focus on the case of $\delta=1$, and the following condition is usually required:
\begin{align}\label{eq:1.6}
\int_0^T \left(\beta_t(\omega)+\gamma^2_t(\omega)\right){\rm d}t<+\infty,\ \ \ps.
\end{align}
The reader is referred to \cite{El KarouiandHuang2011book,BenderKohlmann2000,Yong2006JMAA,
WRC2007JAMSA,LiLaiFan2023PUQR,Li2025,LiFan2026SD,Zhang2026,
Zhu2026,ZhangZhuFan2026} for details. In particular, under the assumptions that the generator $g$ satisfies a stochastic monotonicity condition with a general growth in $y$ and a stochastic Lipschitz continuity condition in $z$ and that the stochastic coefficients $\beta_\cdot$ and $\gamma_\cdot$ only satisfy \eqref{eq:1.6}, \citet{Li2025} and \citet{LiFan2026SD} established existence and uniqueness of the adapted solutions of general time-interval BSDEs in a novel weighted $L^p~(p>1)$ space associated with stochastic coefficients $\beta_\cdot$ and $\gamma_\cdot$. These results were further generalized to the case that the stochastic process $\beta_\cdot$ can take negative values in \citet{Zhang2026}, \citet{Zhu2026} and \citet{ZhangZhuFan2026}.\vspace{0.2cm}

The present paper is devoted to solvability of the unbounded solutions of general time-interval BSDEs when the generator $g$ satisfies \eqref{eq:1.2} with $\delta\in (1,2)$, and the stochastic coefficients $\beta_\cdot\in\R$ and $\gamma_\cdot\in\R_+$ only satisfies the following integrability condition:
\begin{align}\label{eq:1.7}
\int_0^T \left(|\beta_t(\omega)|+\gamma^2_t(\omega)\right){\rm d}t<+\infty,\ \ \ps.
\end{align}
Our results strengthen the corresponding ones of \citet{FH1} and \citet{GWF} to the case that both $\beta_\cdot$ and $\gamma_\cdot$ are possibly unbounded stochastic processes only satisfying \eqref{eq:1.7} and $\beta_\cdot$ can take negative values. First, by constructing a test function $\psi(x,y):=\exp(xy^{2\over \delta^*})$ and applying the It\^{o}-Tanaka formula to $\psi(\nu_t^{1-\frac{2}{\delta^*}(\omega)}, \bar Y_t(\omega))$, where $(\nu_t)_{t\in\T}$ and $(\bar Y_t)_{t\in\T}$ are defined respectively in \eqref{eq:2.2} and \eqref{eq:2.8-1}, we establish a crucial a priori estimate on the unbounded solutions of the preceding BSDEs in a general sub-exponential moment integrability space linked with the stochastic coefficients $\beta_\cdot$ and $\gamma_\cdot$; see \cref{pro:2.2} in \cref{sec:2-Existence} and its proof in \cref{sec:3-Stability theorem} for more details. Based on this estimate and the localization method, we prove existence of the unbounded solution in the preceding space for the case that the generator $g$ satisfies \eqref{eq:1.2} and a general growth condition \ref{A:A2'} along with $\beta_\cdot\geq 0$, and by using an exponential shift transform we address the other case that the generator $g$ satisfies \eqref{eq:1.2} and a slightly stronger growth condition \ref{A:A2} along with $\beta_\cdot\in \R$; see \cref{thm:2.4} in \cref{sec:2-Existence} and its proof in \cref{sec:3-Stability theorem} for more details. Then, we introduce a new general extended convexity/concavity condition related to the sub-quadratic growth of the generator $g$, which generalizes \eqref{eq:1.4} and \eqref{eq:1.5} to the case of unbounded stochastic coefficients; see assumption \ref{A:A3} in \cref{sec:2-Existence} for details. Two propositions (see \cref{pro:3.2,pro:3.3}) are subsequently presented to provide some sufficient conditions such that assumption \ref{A:A3} holds. Finally, by virtue of $\theta$-difference technique we prove a comparison theorem of the unbounded solutions of two BSDEs when either of generators satisfies \ref{A:A3},
yielding naturally the uniqueness of the adapted solution in a general sub-exponential moment integrability space; see \cref{thm:3.5,thm:3.7} in \cref{sec:2-Existence} along with the proof of \cref{thm:3.5} in \cref{sec:3-Stability theorem} for details. Finally, notable corollaries, remarks and examples are included to illustrate our theoretical results; see \cref{cor:2.6,cor:2.8}, \cref{rmk:2.1,rmk:2.3,rmk:2.5,rmk:2.7,rmk:2.9,rmk:3.1,rmk:3.4,rmk:3.6,rmk:3.8,rmk:2.11} and \cref{ex:2.10,ex:3.9} in \cref{sec:2-Existence} along with \cref{rmk:2.5*} in \cref{sec:3-Stability theorem} for more details. We would like to especially emphasize that since the unbounded stochastic coefficient $\beta_\cdot$ can take negative values, the required integrability condition of the terminal value may be weaker than that under the usual situations, which can be indicated through some of these examples.\vspace{0.2cm}

The rest of the paper is organized as follows. In section 2, we state notations and main results, and in section 3 we give the proofs of some theorems and propositions.

\section{Notations and Main results}
\label{sec:2-Existence}
\setcounter{equation}{0}

\subsection{Notations\vspace{0.2cm}}

In this subsection we introduce some notations used later. Let $a\vee b=\max\{a,b\}$, $a^+:=a\vee0$ and $a^-:=(-a)^+$. Let $|x|$ denote the Euclidian norm for $x\in\mathbb{R}^d$ and $x\cdot y$ the usual scalar inner product for $x,y\in\mathbb{R}^d$. Let $\R_+:=[0,+\infty)$, $\R_-:=(-\infty,0]$, $\mathbf{1}_A$ the indicator function of set $A$, and $\mathrm{sgn}(x):= \mathbf{1}_{\{x>0\}}-\mathbf{1}_{\{x\leq 0\}}$. Moreover, by $L^\infty(\Omega,\mathcal{F}_T,\mathbb{P})$ we denote the set of $\mathcal{F}_T$-measurable real-valued bounded random variables, by $\mathcal{S}^\infty$ the set of all $(\F_t)$-progressively measurable and continuous real-valued bounded processes, and by $\mathcal{M}^2$ the set of all $(\F_t)$-progressively measurable $\mathbb{R}^d$-valued processes $(Z_t)_{t\in[0,T]}$ such that
$$
\|Z\|_{\mathcal{M}^2}:=\E\left[\int_0^T|Z_t|^2\mathrm{d}t\right]<+\infty.\vspace{0.2cm}
$$

Throughout the paper, we always suppose that $\delta\in (1,2)$ is a given constant and $\delta^\ast>2$ is the conjugate of $\delta$, i.e., $1/\delta+1/\delta^\ast=1$, and suppose that $(f_t)_{t\in [0,T]}$, $(\gamma_t)_{t\in [0,T]}$ and $(\kappa_t)_{t\in [0,T]}$ are nonnegative real-valued processes and $(\beta_t)_{t\in \T}$ is a real-valued process such that
\begin{equation}\label{eq:2.1}
\int_0^T(f_t(\omega)+|\beta_t(\omega)|+\gamma_t^{\frac{2}{2-\delta}}(\omega))
\, \mathrm{d}t<+\infty
\end{equation}
and
\begin{equation}\label{eq:3.1}
\bigg\|\int_0^T\kappa_t\, \mathrm{d}t\bigg\|_{\infty}\leq \tilde K\vspace{0.2cm}
\end{equation}
with $\tilde K\geq 0$ being a constant. For narrative convenience, we denote for each $t\in\T$,
\begin{align}\label{eq:2.2}
\begin{split}
&\mu_t(\omega):=\int_0^t\beta_s(\omega)\, \mathrm{d}s,\\
&\nu_t(\omega):=C_\delta\int_0^t \left(e^{-(\delta-1)\mu_s(\omega)}\gamma_s(\omega)
\right)^{\frac{2}{2-\delta}}\mathrm{d}s+1,\\
&\bar\nu_t(\omega):=C_\delta\int_0^t \gamma^{\frac{2}{2-\delta}}_s(\omega)
\mathrm{d}s+1
\end{split}
\end{align}
with\vspace{-0.1cm}
\begin{align}\label{eq:2.3}
C_\delta:=(\delta-1)\left(\frac{\delta^2}{\delta-1}\right)^{\frac{\delta}{2-\delta}}.\vspace{0.1cm}
\end{align}
Note that when $\beta_\cdot=0$, $\mu_\cdot=0$ and then $\nu_\cdot$ is identical to $\bar\nu_\cdot$. It follows from \eqref{eq:2.1} and \eqref{eq:2.2} that $(\mu_t)_{t\in\T}$, $(\nu_t)_{t\in\T}$ and $(\bar\nu_t)_{t\in\T}$ are well defined and that for each $t\in\T$,
\begin{align}\label{eq:2.4}
\begin{split}
&\nu_t,\bar\nu_t\geq 1,\\
&\mathrm{d}\nu_t={C}_\delta e^{\frac{2-2\delta}{2-\delta}\mu_t}\gamma_t^{\frac{2}{2-\delta}}\mathrm{d}t,\\
&\mathrm{d}\bar\nu_t={C}_\delta\gamma_t^{\frac{2}{2-\delta}}\mathrm{d}t.
\end{split}
\end{align}
Throughout the paper, we also always suppose that $\xi$ is a real-valued terminal condition and $g$ is a generator of BSDE \eqref{eq:1.1} which is continuous in $(y,z)$. \vspace{0.1cm}

\subsection{A priori estimate\vspace{0.2cm}}

Let us first introduce the following assumption of the generator $g$.

\begin{enumerate}
\renewcommand{\theenumi}{(A1)}
\renewcommand{\labelenumi}{\theenumi}
\item\label{A:A1} $g$ has a stochastic one-sided linear growth in $y$ and a stochastic sub-quadratic growth in $z$, i.e., $\as$, for each $(y,z)\in \R\times\R^d$, we have
$$\mathrm{sgn}(y)g(\omega,t,y,z)\leq f_t(\omega)+\beta_{t}(\omega)|y|+\gamma_t(\omega)|z|^\delta.$$
\end{enumerate}

\cref{pro:2.2} below establishes an a priori estimate for adapted solutions of BSDEs with sub-quadratic generators and unbounded stochastic coefficients, which will play a crucial role in the proof of the main result of the paper. Its proof will be given in next section.

\begin{pro}\label{pro:2.2}
Assume that $\xi$ is a terminal condition, $g$ is a generator satisfying \ref{A:A1}, and $(Y_t,Z_t)_{t\in[0,T]}$ is an adapted solution to BSDE $(\xi,T,g)$ such that
\begin{equation}\label{eq:2.5}
\E\left[\sup_{t\in\T}\exp\left(\nu_t^{1-\frac{2}{\delta^\ast}}\left(e^{\mu_t}
|Y_t|+(\delta^\ast)^{\frac{\delta^\ast}{2}}+\int_0^t e^{\mu_s}f_s\mathrm{d}s
\right)^{\frac{2}{\delta^\ast}}\right)\right]<+\infty.
\end{equation}
Then, we have for each $t\in[0,T]$,
\begin{align}\label{eq:2.6}
\exp{\left(\nu_t^{1-\frac{2}{\delta^\ast}}(e^{\mu_t}|Y_t|)^{\frac{2}{\delta^\ast}}\right)}
\leq\mathbb{E}\left[\exp\left(\nu_T^{1-\frac{2}{\delta^\ast}}
\left(e^{\mu_T}|\xi|+(\delta^\ast)^{\frac{\delta^\ast}{2}}+
\displaystyle\int_0^Te^{\mu_s}f_s\mathrm{d}s\right)^{\frac{2}{\delta^\ast}}\right)
\bigg|\mathcal{F}_t\right].
\end{align}
\end{pro}\vspace{0.02cm}

\begin{rmk}\label{rmk:2.3}
According to the proof of \cref{pro:2.2}, we can verify the following assertion: in \cref{pro:2.2}, if $|\xi|$ is replaced with $\xi^+$, the assumption \ref{A:A1} is replaced with
\begin{enumerate}
\renewcommand{\theenumi}{(A1')}
\renewcommand{\labelenumi}{\theenumi}
\item\label{A:A1'} \ \ $\as,\ \
\mathbf{1}_{\{y>0\}}g(\omega,t,y,z)\leq f_t(\omega)+\beta_{t}(\omega)y^+ +\gamma_t(\omega)|z|^\delta,~~(y,z)\in \mathbb{R}\times\mathbb{R}^d,$
\end{enumerate}
and the condition \eqref{eq:2.5} is replaced with
$$
\E\left[\sup_{t\in\T}\exp\left(\nu_t^{1-\frac{2}{\delta^\ast}}\left(e^{\mu_t}
Y_t^+ +(\delta^\ast)^{\frac{\delta^\ast}{2}}+\int_0^t e^{\mu_s}f_s\mathrm{d}s
\right)^{\frac{2}{\delta^\ast}}\right)\right]<+\infty,
$$
then for each $t\in[0,T]$, we have
$$\exp\left(\nu_t^{1-\frac{2}{\delta^\ast}}(e^{\mu_t}Y_t^+)^{\frac{2}{\delta^\ast}}\right)\leq\mathbb{E}\left[
\exp\left(\nu_T^{1-\frac{2}{\delta^\ast}}\left(e^{\mu_T}\xi^+
+(\delta^\ast)^{\frac{\delta^\ast}{2}}+\int_0^Te^{\mu_s}f_s\mathrm{d}s\right)
^{\frac{2}{\delta^\ast}}\right)\bigg|\mathcal{F}_t\right].$$
Indeed, we only need to substitute $|Y_\cdot|,~\mathrm{sgn}(Y_\cdot)$ and $L_\cdot$ in the proof of \cref{pro:2.2} into $Y_\cdot^+, \mathbf{1}_{\{Y_\cdot>0\}}$ and $L_\cdot/2$ respectively. \vspace{0.1cm}

We would like to mention that \cref{pro:2.2} and the preceding assertion generalize the corresponding conclusions stated in \citet{FH1} and \citet{GWF}, where both $\beta_\cdot$ and $\gamma_\cdot$ are either nonnegative constants or nonnegative deterministic functions.\vspace{0.1cm}

We also would like to mention that for the case of $\delta=2$, the proof of \cref{pro:2.2} is not valid any longer. For this case, in order to obtain the desired a priori estimate, it seems that the process $\gamma_\cdot$ needs to be bounded by a nonnegative constant.
\vspace{0.2cm}
\end{rmk}

\subsection{Existence\vspace{0.2cm}}

Let us further introduce the following assumption of the generator $g$.

\begin{enumerate}
\renewcommand{\theenumi}{(A2)}
\renewcommand{\labelenumi}{\theenumi}
\item\label{A:A2} $g$ has a general growth in $y$ and a sub-quadratic growth in $z$, i.e., there exists  a function $\varphi_t(\omega,x)\in \bs$ such that $\as$, for each $(y,z)\in \R\times\R^d$, we have
$$
|g(\omega,t,y,z)|\leq \varphi_t(\omega,|y|)+\gamma_t(\omega)|z|^{\delta}.
$$
here and hereafter, $\bs$ denotes the set of all $\varphi_t(\omega,x):\Omega\times[0,T]\times \R_+\rightarrow\R_+$ satisfying
\begin{enumerate}
\item [(i)] Almost everywhere in $\Omega\times\T$, $\varphi_t(\omega,\cdot)$ is increasing;
\item [(ii)] For each $r>0$, we have
$$\int_0^T\varphi_t(\omega,r)\mathrm{d}t<+\infty.$$
\end{enumerate}

\renewcommand{\theenumi}{(A2')}
\renewcommand{\labelenumi}{\theenumi}
\item\label{A:A2'} $g$ has a general growth in $y$ and a quadratic growth in $z$, i.e., there exist an increasing function $\phi(x):\R_+\To \R_+$ and a function $\bar\varphi_t(\omega,x)\in \bs$ such that $\as$, for each $(y,z)\in \R\times\R^d$, we have
$$
|g(\omega,t,y,z)|\leq \bar\varphi_t(\omega,|y|)+\phi(|y|)|z|^2.
$$
\end{enumerate}

\begin{rmk}\label{rmk:2.1}
It follows from Young's inequality that for each $(\omega,t,z)\in\Omega\times\T\times\R^d$,
$$
\gamma_t|z|^\delta\leq \frac{2-\delta}{2} \gamma_t^{\frac{2}{2-\delta}}+\frac{\delta}{2}|z|^2.
$$
Thus, in light of \eqref{eq:2.1}, if the generator $g$ satisfies (A2), then it also satisfies (A2') with $\phi(x):=\frac{\delta}{2}$ and $\bar\varphi_t(x):=\varphi_t(x)+\frac{2-\delta}{2} \gamma_t^{\frac{2}{2-\delta}}\in\bs$. That is to say, (A2') is weaker than (A2).\vspace{0.2cm}
\end{rmk}

Now, we give the following existence theorem which is one of the main results of the paper. Its proof will be provided in next section.

\begin{thm}\label{thm:2.4}
Assume that $\xi$ is a terminal condition and $g$ is a generator satisfying assumptions \ref{A:A1} and \ref{A:A2}. If the following two conditions
\begin{equation}\label{eq:2.18}
\E\left[\exp\left(k\nu_T^{1-\frac{2}{\delta^\ast}}\right)\right]<+\infty,\ \ \RE k>0,
\end{equation}
and
\begin{equation}\label{eq:2.19}
\mathbb{E}\left[\exp\left(\nu_T^{1-\frac{2}{\delta^\ast}}\left(e^{\mu_T}
|\xi|+(\delta^\ast)^{\frac{\delta^\ast}{2}}+\displaystyle\int_0^Te^{\mu_s}f_s\mathrm{d}s
\right)^{\frac{2}{\delta^\ast}}\right)\right]<+\infty\vspace{0.2cm}
\end{equation}
are fulfilled, then BSDE $(\xi,T,g)$ admits a solution $(Y_t,Z_t)_{t\in[0,T]}$ such that for each $t\in[0,T]$,
\begin{align}\label{eq:2.6*}
\exp\left(\nu_t^{1-\frac{2}{\delta^\ast}}(e^{\mu_t}|Y_t|)^{\frac{2}{\delta^\ast}}\right)
\leq\mathbb{E}\left[\exp\left(\nu_T^{1-\frac{2}{\delta^\ast}}
\left(e^{\mu_T}|\xi|+(\delta^\ast)^{\frac{\delta^\ast}{2}}+
\displaystyle\int_0^Te^{\mu_s}f_s\mathrm{d}s\right)^{\frac{2}{\delta^\ast}}\right)
\bigg|\mathcal{F}_t\right].
\end{align}

Moreover, if $(\beta_t)_{t\in\T}$ is a nonnegative process, then the above assertion remains true when the assumption \ref{A:A2} is replaced with the weaker assumption \ref{A:A2'}.\vspace{0.2cm}
\end{thm}

\begin{rmk}\label{rmk:2.5}
It is evident that if the generator $g$ has a stochastic linear growth in $y$ and a stochastic sub-quadratic growth in $z$, i.e., $\as$, for each $(y,z)\in \R\times\R^d$, we have
\begin{align*}
|g(\omega,t,y,z)|\leq f_t(\omega)+\beta_{t}(\omega)|y|+\gamma_t(\omega)|z|^\delta,
\end{align*}
then $g$ must satisfy \ref{A:A1} and \ref{A:A2}. By the way, $\beta_\cdot$ has to be nonnegative in this case.\vspace{0.2cm}
\end{rmk}

The following \cref{cor:2.6} is a direct consequence of \cref{thm:2.4}.

\begin{cor}\label{cor:2.6}
Assume that $\xi$ is a terminal condition, $g$ is a generator satisfying assumptions \ref{A:A1} and \ref{A:A2}, and that for some constant $M>0$,\vspace{0.1cm}
\begin{equation}\label{eq:2.26}
 \int_0^T(|\beta_s|+\gamma_s^{\frac{2}{2-\delta}})\mathrm{d}s\leq M.\vspace{-0.2cm}
\end{equation}
If $|\xi|+\displaystyle\int_0^Tf_s\mathrm{d}s$ has sub-exponential moments of any order, i.e.,
\begin{equation}\label{eq:2.27}
\mathbb{E}\left[\exp\left(p\left(|\xi|+\displaystyle\int_0^Tf_s\mathrm{d}s
\right)^{\frac{2}{\delta^\ast}}\right)\right]<+\infty,~~\forall p>1
\end{equation}
then BSDE $(\xi,T,g)$ admits an adapted solution $(Y_t,Z_t)_{t\in[0,T]}$ such that
\begin{equation}\label{eq:2.28}
\exp\left(|Y_t|^{\frac{2}{\delta^\ast}}\right)
\leq K\mathbb{E}\left[\left.\exp\left(K\left(|\xi|+\displaystyle\int_0^Tf_s\mathrm{d}s
\right)^{\frac{2}{\delta^\ast}}\right)\right|\F_t\right],~~t\in[0,T],
\end{equation}
where $K>0$ is a constant only depending on $(T,\delta,M)$.\vspace{0.2cm}
\end{cor}

\begin{proof}
By \eqref{eq:2.26} and \eqref{eq:2.27}, the definition of $\mu_\cdot$ and $\nu_\cdot$, and the inequality $(a+b)^\lambda\leq a^\lambda+b^\lambda$ for each $a,b>0$ and $\lambda\in(0,1)$, we have
$$
\E\left[\exp\left(k\nu_T^{1-\frac{2}{\delta^\ast}}\right)\right]\leq \exp\left(k\left(TC_\delta Me^{\frac{2\delta-2}{2-\delta}M}+1\right)^{1-\frac{2}{\delta^\ast}}\right)=:\exp\left(k\bar C_{T,\delta,M}\right),\ \ \RE k>0
$$
and
\begin{align}\label{eq:2.29}
\begin{split}
\exp&\left(\nu_T^{1-\frac{2}{\delta^\ast}}\left(e^{\mu_T}|\xi|+(\delta^\ast)^{\frac{\delta^\ast}{2}}
+\displaystyle\int_0^Te^{\mu_s}f_s\mathrm{d}s\right)^{\frac{2}{\delta^\ast}}\right)\\
&\leq\exp\left(\bar C_{T,\delta,M}\left(\delta^\ast
+\exp\left({\frac{2}{\delta^\ast}M}\right)\left(|\xi|+\int_0^Tf_s\mathrm{d}s\right)^{\frac{2}{\delta^\ast}}\right)\right)\\
&\leq K\exp\left(K\left(|\xi|+\int_0^Tf_s\mathrm{d}s\right)^{\frac{2}{\delta^\ast}}\right),
\end{split}
\end{align}
where
$$
K:=\exp\left(\delta^\ast\bar C_{T,\delta,M}\right)\vee\left(\bar C_{T,\delta,M}\exp\left({\frac{2}{\delta^\ast}M}\right)\right).\vspace{0.2cm}
$$
Thus, we have \eqref{eq:2.18} and \eqref{eq:2.19}. It then follows from \cref{thm:2.4} that BSDE $(\xi,T,g)$ admits a solution $(Y_t,Z_t)_{t\in[0,T]}$ such that $Y_\cdot$ satisfies \eqref{eq:2.6*}, and \eqref{eq:2.28} follows by \eqref{eq:2.6*} and \eqref{eq:2.29}. \vspace{0.2cm}
\end{proof}

\begin{rmk}\label{rmk:2.7}
By a similar argument to \eqref{eq:2.24} in the proof of \cref{thm:2.4}, it is not difficult to verify that for the case of $\int_0^T |\beta_t|{\rm d}t\leq M$ for some $M>0$, if the generator $g$ satisfies \ref{A:A2'}, then the generator $\tilde g$ defined via \eqref{eq:2.24*} also satisfies \ref{A:A2'} with a different $\phi(\cdot)$.  Hence, by the proof of \cref{thm:2.4} we can conclude that the assertion of \cref{cor:2.6} remains true when \ref{A:A2} is weakened to \ref{A:A2'}. This strengthens \citet[Theorem 3.2]{GWF}, where both $\beta_\cdot$ and $\gamma_\cdot$ are nonnegative deterministic functions.\vspace{0.2cm}
\end{rmk}

\cref{cor:2.8} below is another direct consequence of \cref{thm:2.4}. We omit its proof here.

\begin{cor}\label{cor:2.8}
Assume that $\xi$ is a terminal condition, $g$ is a generator satisfying assumptions \ref{A:A1} and \ref{A:A2}, and $\nu_T\leq M$ for some constant $M>0$. If
\begin{equation*}
\mathbb{E}\left[\exp\left(p\left(e^{\mu_T}|\xi|+\displaystyle\int_0^T e^{\mu_s}f_s\mathrm{d}s
\right)^{\frac{2}{\delta^\ast}}\right)\right]<+\infty,~~\forall p>1
\end{equation*}
then BSDE $(\xi,T,g)$ admits an adapted solution $(Y_t,Z_t)_{t\in[0,T]}$ such that
\begin{equation*}
\exp(({e^{\mu_t}|Y_t|)^{\frac{2}{\delta^\ast}})}
\leq K\mathbb{E}\left[\left.\exp\left(K\left(e^{\mu_T}|\xi|+\int_0^T e^{\mu_s}f_s\mathrm{d}s
\right)^{\frac{2}{\delta^\ast}}\right)\right|\F_t\right],~~t\in[0,T],
\end{equation*}
where $K>0$ is a constant only depending on $\delta$ and $M$.\vspace{0.2cm}
\end{cor}

\begin{rmk}\label{rmk:2.9}
Suppose that there exists a deterministic function $h(t)\in \T\To \R_+$ satisfying $\int_0^T h^{\frac{2}{2-\delta}}(t){\rm d}t<+\infty$ such that
$$
e^{-(\delta-1)\mu_t}\gamma_t\leq h(t),\ \ t\in\T.
$$
Then, by the definition of $\nu_\cdot$ we have that $\nu_T\leq M$ for some constant $M>0$.\vspace{0.2cm}
\end{rmk}

\begin{ex}\label{ex:2.10}
We provide several examples to demonstrate that \cref{thm:2.4}, \cref{cor:2.6} and  \cref{cor:2.8} are all not encompassed by any existing conclusions.
\begin{itemize}
\item [(i)] For each $(\omega,t,y,z)\in \Omega\times\T\times \R\times\R^{d}$, let
$$
g(\omega,t,y,z):=\left(|B_t(\omega)|^{1\over 4}{\bf 1}_{\{0\leq t\leq 1\}}+e^{-t}\right)|z|^{3\over 2}\sin|z|+|z|^2 (1-e^y).
$$
It is not very difficult to verify that this $g$ satisfies \ref{A:A1} and \ref{A:A2'} with
$$
f_t=\beta_t:=0,\ \ \gamma_t:=|B_t|^{1\over 4}{\bf 1}_{\{0\leq t\leq 1\}}+e^{-t},\ \ \delta:=\frac{3}{2},\ \ \bar\varphi_t(x):=\gamma_t^4\ \ {\rm and}\ \ \phi(x):=e^{x}+2.
$$
Note that
$$\nu_T=\int_0^T \gamma^4_t\, {\rm d}t\leq 8\int_0^1|B_t|\, {\rm d}t+2,$$
and then $\E[\exp(k\nu_T^{1\over 3})]<+\infty$ for each $k>0$. It then follows from \cref{thm:2.4} that for each $\mathcal{F}_T$-measurable random variable $\xi$ satisfying that
$$\E\left[\exp\left(\left(8\int_0^1|B_t|\, {\rm d}t+2\right)^{1\over 3}|\xi|^{2\over 3}\right)\right]<+\infty,$$
BSDE $(\xi,T,g)$ admits an adapted solution $(Y_t,Z_t)_{t\in \T}$ such that
$$
\exp\left(|Y_t|^{2\over 3}\right)\leq \E\left[\left.\exp\left(\left(8\int_0^1|B_t|\, {\rm d}t+2\right)^{1\over 3}\left(|\xi|+\sqrt{27}\,\right)^{2\over 3}\right)\right|\F_t\right],\ \ t\in\T.
$$

\item [(ii)] Let $T=+\infty$ and $M>0$ be a constant. For each $(\omega,t,y,z)\in \Omega\times\T\times \R\times\R^{d}$, let
$$
g(\omega,t,y,z):=|B_t|{\bf 1}_{\{0\leq t\leq \sigma_1\}}(|y|+1)+|B_t|^2{\bf 1}_{\{0\leq t\leq \sigma_2\}}|z|^{6\over 5}\sin y-y^3|z|^2,
$$
where $\sigma_1$ and $\sigma_2$ are two $(\F_t)$-stopping times defined by
$$
\sigma_1:=\inf\left\{t\geq 0: \int_0^t |B_s|\, {\rm d}s\geq M\right\}
$$
and
$$
\sigma_2:=\inf\left\{t\geq 0: \int_0^t |B_s|^5\, {\rm d}s\geq M\right\}.\vspace{0.2cm}
$$
It is not very difficult to verify that this $g$ satisfies \ref{A:A1} and \ref{A:A2'} with
$$
f_t=\beta_t:=|B_t|{\bf 1}_{\{0\leq t\leq \sigma_1\}},\ \ \gamma_t:=|B_t|^2{\bf 1}_{\{0\leq t\leq \sigma_2\}},\ \ \delta:=\frac{6}{5},\ \ \bar\varphi_t(x):=(1+x)f_t+\gamma_t^{5\over 2}
$$
and $\phi(x):=x^3+1$. Note by the definitions of $\sigma_1$ and $\sigma_2$ that $\int_0^T (|\beta_t|+\gamma^{5\over 2}_t)\, {\rm d}t\leq 2M$. It then follows from \cref{cor:2.6} and \cref{rmk:2.7} that for each $\mathcal{F}_T$-measurable random variable $\xi$ satisfying that
$$\E\left[\exp\left(k|\xi|^{1\over 3}\right)\right]<+\infty,\ \ \RE k>0,$$
BSDE $(\xi,T,g)$ admits an adapted solution $(Y_t,Z_t)_{t\in \T}$ such that
$$
\exp\left(|Y_t|^{1\over 3}\right)\leq K\E\left[\left.\exp\left(K|\xi|^{1\over 3}\right)\right|\F_t\right],\ \ t\in\T,
$$
where $K>0$ is a constant depending only on $M$.

\item [(iii)] Let $T$ be a finite positive real number, and for each $(\omega,t,y,z)\in \Omega\times\T\times \R\times\R^{d}$, let
$$
g(\omega,t,y,z):=-y^3 e^{|B_t|^2}+|B_t| y\cos y+e^{\frac{1}{3}\int_0^t |B_s|{\rm d}s}|z|^{4\over 3}\sin(y^2|z|).
$$
It is not very difficult to verify that this $g$ satisfies \ref{A:A1} and \ref{A:A2} with
$$
f_t:=0,\ \ \beta_t:=|B_t|,\ \ \gamma_t:=e^{\frac{1}{3}\int_0^t |B_s|{\rm d}s},\ \ \delta:=\frac{4}{3}\ \ {\rm and}\ \ \varphi_t(x):=e^{|B_t|^2} x^3+|B_t|x.
$$
Note that $e^{-(\delta-1)\mu_t}\gamma_t\leq 1$ for each $t\in\T$ and then $\nu_T\leq T$. It then follows from \cref{cor:2.8} that for each $\mathcal{F}_T$-measurable random variable $\xi$ satisfying that
$$\E\left[\exp\left(p\left(e^{\int_0^T|B_t|\,{\rm d}t}|\xi|\right)^{1\over 2}\right)\right]<+\infty,\ \ \RE p>0,$$
BSDE $(\xi,T,g)$ admits an adapted solution $(Y_t,Z_t)_{t\in \T}$ such that
$$
\exp\left(\left(e^{\int_0^t|B_s|\,{\rm d}s}|Y_t|\right)^{1\over 2}\right)\leq K\E\left[\left.\exp\left(K\left(e^{\int_0^T|B_t|\,{\rm d}t}|\xi|\right)^{1\over 2}\right)\right|\F_t\right],\ \ t\in\T,
$$
where $K>0$ is a constant depending only on $T$.

\item [(iv)] Let $T$ be a finite positive real number, and for each $(\omega,t,y,z)\in \Omega\times\T\times \R\times\R^{d}$, let
$$
g(\omega,t,y,z):=-|B_t|^3 y +2e^{-\frac{1}{4}\int_0^t |B_s|^3{\rm d}s}|z|^{5\over 4}\cos |z|.
$$
It is not very difficult to verify that this $g$ satisfies \ref{A:A1} and \ref{A:A2} with
$$
f_t:=0,\ \ \beta_t:=-|B_t|^3,\ \ \gamma_t:=2e^{-\frac{1}{4}\int_0^t |B_s|^3{\rm d}s},\ \ \delta:=\frac{5}{4}\ \ {\rm and}\ \ \varphi_t(x):=|B_t|^3 x.
$$
Note that $e^{-(\delta-1)\mu_t}\gamma_t\leq 2$ for each $t\in\T$ and then $\nu_T\leq 2T$. It follows from \cref{cor:2.8} that for $\xi:=e^{\int_0^T |B_t|^3\, {\rm d}t}$, BSDE $(\xi,T,g)$ admits a solution $(Y_t,Z_t)_{t\in \T}$ such that
$$
|Y_t|\leq K e^{\int_0^t |B_s|^3\, {\rm d}s},\ \ t\in\T,\vspace{-0.2cm}
$$
where $K>0$ is a constant depending only on $T$. We would like to especially mention that this $\xi$ is not an integrable random variable.\vspace{0.2cm}
\end{itemize}
\end{ex}

\begin{rmk}\label{rmk:2.11}
We note that the integrability assertion of $(Z_t)_{t\in\T}$ presented in \citet[Theorem 3.2]{GWF} is not obtained in \cref{thm:2.4}, \cref{cor:2.6,cor:2.8}. The main reason lies in the less integrability requirement of unbounded stochastic coefficients $\beta_\cdot$ and $\gamma_\cdot$.
\end{rmk}

\subsection{Comparison theorem\vspace{0.2cm}}

In this subsection, we will establish a comparison theorem for unbounded solutions of sub-quadratic BSDEs with unbounded stochastic coefficients under a general extended convexity/concavity condition related to the stochastic sub-quadratic growth. \vspace{0.2cm}

Let us first introduce the following general extended convexity/concavity assumption of the generator $g$ related to the stochastic sub-quadratic growth.

\begin{enumerate}
\renewcommand{\theenumi}{(A3)}
\renewcommand{\labelenumi}{\theenumi}
\item\label{A:A3} $g$ satisfies a general extended convexity/concavity condition related to the sub-quadratic growth, i.e., $\as$, for each $(y_i,z_i)\in \R\times\R^d,\, i=1,2$ and each $\theta\in(0,1)$, we have
\begin{align}\label{eq:3.2}
\begin{split}
&\mathbf{1}_{\{\Delta_\theta y>0\}}\frac{g(\omega,t,y_1,z_1)-\theta g(\omega,t,y_2,z_2)}{1-\theta}\\
&\ \ \leq f_t(\omega)+\kappa_t(\omega)|y_2|+\beta_t(\omega)(\Delta_\theta y)^+ +\gamma_t(\omega)|\Delta_\theta z|^\delta
\end{split}
\end{align}
or
\begin{align}\label{eq:3.3}
\begin{split}
&-\mathbf{1}_{\{\Delta_\theta y<0\}}\frac{g(\omega,t,y_1,z_1)-\theta g(\omega,t,y_2,z_2)}{1-\theta}\\
&\ \ \leq f_t(\omega)+\kappa_t(\omega)|y_2|+\beta_t(\omega)(\Delta_\theta y)^- +\gamma_t(\omega)|\Delta_\theta z|^\delta,
\end{split}
\end{align}
where
$$
\Delta_\theta y:=\frac{y_1-\theta y_2}{1-\theta}\ \ {\rm and}\ \ \Delta_\theta z:=\frac{z_1-\theta z_2}{1-\theta}.\vspace{0.2cm}
$$
\end{enumerate}

\begin{rmk}\label{rmk:3.1} With respect to assumption \ref{A:A3}, we make the following comments.
\begin{itemize}
\item [(i)] Similar assumptions to \ref{A:A3} are first presented in \citet{FanHuTang2020CRM} and then in \citet{FH1, GWF,FanHuTang2026MF}, where $\beta_\cdot$, $\gamma_\cdot$ and $\kappa_\cdot$ are all either nonnegative constants or deterministic nonnegative functions.

\item [(ii)] Letting $y_1=y_2=y$ and $z_1=z_2=z$ in \eqref{eq:3.2} and \eqref{eq:3.3} respectively yields that
    $$
    {\bf 1}_{\{y>0\}} g(\omega,t,y,z)\leq f_t(\omega)+(\kappa_t(\omega)+\beta_t(\omega))y^+ +\gamma_t(\omega)|z|^\delta
    $$
and
    $$
    -{\bf 1}_{\{y<0\}} g(\omega,t,y,z)\leq f_t(\omega)+(\kappa_t(\omega)+\beta_t(\omega))y^- +\gamma_t(\omega)|z|^\delta,
    $$
whose combination implies that $g$ satisfies \ref{A:A1} with $\beta_t$ being replaced with
$\kappa_t+\beta_t$.

\item [(iii)] Letting first $z_1=z_2=z$ in \eqref{eq:3.2} or \eqref{eq:3.3} and then letting $\theta\To 1$ implies that
    $${\rm sgn}(y_1-y_2)(g(\omega,t,y_1,z)-g(\omega,t,y_2,z))\leq \beta_t(\omega) |y_1-y_2|,$$
which indicates that $g$ satisfies the so-called stochastic monotonicity condition in $y$.\vspace{0.2cm}
\end{itemize}
\end{rmk}

By a similar argument to that of \citet[Proposition 3.5]{FH1}, we can prove the following proposition. Its proof is omitted here.

\begin{pro}\label{pro:3.2}
Assume that $\as$, the generator $g$ is continuous in $(y,z)$.
\begin{itemize}
\item [(i)] If $g$ satisfies assumption \ref{A:A1} and $\as$, $g(\omega,t,\cdot,\cdot)$ is convex (resp. concave), then it also satisfies \eqref{eq:3.2} (resp. \eqref{eq:3.3}) with $\kappa_\cdot=0$.

\item [(ii)] If $g$ satisfies assumption \ref{A:A1} and $\as$, for each $(y,z)\in \R\times\R^d$,  $g(\omega,t,\cdot,z)$ is Lipschitz continuous with the Lipschitz coefficient $\beta_t(\omega)\geq 0$ and $g(\omega,t,\cdot,\cdot)$ is convex (resp. concave), then $g$ also satisfies assumption \eqref{eq:3.2} (resp. \eqref{eq:3.3}) with $\kappa_\cdot=2\beta_\cdot$.

\item [(iii)]Suppose that $l(\omega,t,y):\Omega\times\T\times\R\To \R$ and $q(\omega,t,z):\Omega\times\T\times\R^d\To \R$ are $(\F_t)$-progressively measurable processes for each $(y,z)$, and that $\as$, $|l(\omega,t,\cdot)|$ and $|q(\omega,t,\cdot)|$ are bounded by some uniform constant $M>0$ and Lipschitz continuous with the Lipschitz coefficients $\beta_t(\omega)\geq 0$ and $\gamma_t(\omega)\geq 0$, respectively. Define
    $$
    g(\omega,t,y,z):=l(\omega,t,y)q(\omega,t,z),\ \ (\omega,t,y,z)\in \Omega\times\T\times\R\times\R^d.
    $$
    If $\as$, $q(\omega,t,z)=0$ when $|z|\geq R$ for some uniform constant $R>0$, then $\as$, for each $(y_i,z_i)\in \R\times\R^d,\, i=1,2$ and each $\theta\in(0,1)$, we have
\begin{align*}
\begin{split}
&\left|\frac{g(\omega,t,y_1,z_1)-\theta g(\omega,t,y_2,z_2)}{1-\theta}\right|\\
&\ \ \leq M\left(5M+(2R+1)\gamma_t(\omega)+\beta_t(\omega)|y_2|+\beta_t(\omega)|\Delta_\theta y| +\gamma_t(\omega)|\Delta_\theta z|^\delta\right).
\end{split}
\end{align*}
\end{itemize}
\end{pro}\vspace{0.2cm}

\cref{pro:3.3} below is closely connected with the assumption \ref{A:A3}. Its proof will be provided in next section. Some examples of generators satisfying \ref{A:A1}, \ref{A:A2}/\ref{A:A2'} and \ref{A:A3} will be constructed later via \cref{pro:3.2,pro:3.3}.

\begin{pro}\label{pro:3.3}
We have the following assertions.
\begin{itemize}
\item [(i)] Clearly, the statement that the generator $g$ satisfies \eqref{eq:3.3} is equivalent to the statement that the generator $\bar g(\omega,t,y,z):=-g(\omega,t,-y,-z)$ satisfies \eqref{eq:3.2}.

\item [(ii)] Let $g(\omega,t,y):=\beta_t(\omega)y,\ (\omega,t,y)\in\Omega\times\T\times\R$. Then, $\as$, for each $(y_1,y_2)\in \R\times\R$ and each $\theta\in(0,1)$, we have
\begin{align*}
{\bf 1}_{\{\Delta_\theta y>0\}}\frac{g(\omega,t,y_1)-\theta g(\omega,t,y_2)}{1-\theta}=\beta_t(\omega)(\Delta_\theta y)^+.
\end{align*}

\item [(iii)] Assume that $g(\omega,t,y):\Omega\times\T\times\R\To\R$ is an $(\F_t)$-progressively measurable process for each $y\in\R$, and that $\as$, $g(\omega, t,\cdot)$ is non-increasing and continuous on $\R$ and $g(\omega,t,0)=0$. If $\as$, $g(\omega, t,\cdot)$ is both convex on $\R_-$ and $\R_+$, then $\as$, for each $(y_1,y_2)\in \R\times\R$ and each $\theta\in(0,1)$, we have
\begin{align}\label{eq:3.4}
{\bf 1}_{\{\Delta_\theta y>0\}}\frac{g(\omega,t,y_1)-\theta g(\omega,t,y_2)}{1-\theta}\leq 0;
\end{align}
    If $\as$, $g(\omega, t,y)\leq \kappa_t(\omega)|y|$ for each $y\in\R_-$ and $g(\omega, t,\cdot)$ is convex on $\R_+$, then $\as$, for each $(y_1,y_2)\in \R\times\R$ and each $\theta\in(0,1)$, we have
\begin{align}\label{eq:3.5}
{\bf 1}_{\{\Delta_\theta y>0\}}\frac{g(\omega,t,y_1)-\theta g(\omega,t,y_2)}{1-\theta}\leq \kappa_t(\omega)|y_2|.
\end{align}

\item [(iv)] Assume that $a>0$ and $b\in\R$ are two given constants, and that the generator $g(\omega,t,y):\Omega\times\T\times\R\To \R$ satisfies \eqref{eq:3.2}. Define
    $$
    \bar g(\omega,t,y):=g(\omega,t,y+a)+b,\ \ (\omega,t,y)\in \Omega\times\T\times\R.
    $$
    Then $\as$, for each $(y_1,y_2)\in \R\times\R$ and each $\theta\in(0,1)$, we have
\begin{align*}
\begin{split}
&{\bf 1}_{\{\Delta_\theta y>0\}}\frac{\bar g(\omega,t,y_1)-\theta \bar g(\omega,t,y_2)}{1-\theta}\\
&\ \ \leq f_t(\omega)+b+(\kappa_t(\omega)+|\beta_t(\omega)|)a+\kappa_t(\omega)|y_2|+\beta_t(\omega)(\Delta_\theta y)^+.
\end{split}
\end{align*}

\item [(v)] It is not difficult to verify that if both $g_1$ and $g_2$ satisfy \eqref{eq:3.2} (resp. \eqref{eq:3.3}), then for each $k_1,k_2\geq 0$, the generator $k_1g_1+k_2g_2$ also satisfies \eqref{eq:3.2} (resp. \eqref{eq:3.3}), but with different parameters.\vspace{0.2cm}
\end{itemize}
\end{pro}

\begin{rmk}\label{rmk:3.4}
Concerning \cref{pro:3.3}, we have the following comments.
\begin{itemize}
\item [(i)] For each $(\omega,t,y)\in \Omega\times\T\times\R$, define
$$
g(\omega,t,y):=y^2{\bf 1}_{\{y\leq 0\}}-\sqrt{y}{\bf 1}_{\{y>0\}}\ \ {\rm and}\ \ \bar g(\omega,t,y):=-\kappa_t(\omega)\sqrt[3]{y}.
$$
According to (iii) of \cref{pro:3.3}, we can conclude that $g$ satisfies \eqref{eq:3.4} and $\bar g$ satisfies \eqref{eq:3.5}. This indicates that $g$ satisfying \ref{A:A3} can have a general growth, be non-convex/non-concave and not locally Lipschitz continuous in $y$. In addition, it follows from (iii) of \cref{pro:3.2} that $g$ satisfying \ref{A:A3} can be non-convex/non-concave in $z$.

\item [(ii)] It should be noted that the assertion in (iv) of \cref{pro:3.3} does not hold in general when $a<0$. In fact, define
    $$
    g(\omega,t,y):=\sqrt{-y}\, {\bf 1}_{\{y\leq 0\}},\ \ (\omega,t,y)\in \Omega\times\T\times\R.
    $$
    It follows from (iii) of \cref{pro:3.3} that $g$ satisfies \eqref{eq:3.2}. More specifically, $\as$, for each $(y_1,y_2)\in \R\times\R$ and each $\theta\in(0,1)$, we have
\begin{align*}
{\bf 1}_{\{\Delta_\theta y>0\}}\frac{g(\omega,t,y_1)-\theta g(\omega,t,y_2)}{1-\theta}\leq |y_2|.
\end{align*}
However, the following generator
$$
\bar g(\omega,t,y):=g(\omega,t,y-1)=\sqrt{-y+1}\, {\bf 1}_{\{y\leq 1\}},\ \ (\omega,t,y)\in \Omega\times\T\times\R
$$
does not satisfy \eqref{eq:3.2}. Indeed, if there exists a constant $c>0$ such that $\as$,
for each $(y_1,y_2)\in \R\times\R$ and each $\theta\in(0,1)$,
\begin{align}\label{eq:3.6}
{\bf 1}_{\{\Delta_\theta y>0\}}\frac{\bar g(\omega,t,y_1)-\theta \bar g(\omega,t,y_2)}{1-\theta}\leq c\left(1+|y_2|+(\Delta_\theta y)^+\right),
\end{align}
then by picking $y_1:=\theta$ and $y_2:=1-\eps$ for some constant $\eps\in (0,1)$, we have
$$
\Delta_\theta y=\frac{\theta-\theta(1-\eps)}{1-\theta}=\frac{\eps\theta}{1-\theta}>0
$$
and
$$
\frac{\bar g(\omega,t,\theta)-\theta \bar g(\omega,t,1-\eps)}{1-\theta}=\frac{\sqrt{-\theta+1}-\theta \sqrt{\eps}}{1-\theta}.
$$
Substituting these two equalities into \eqref{eq:3.6} and letting $\eps\To 0$ yields that for each $\theta\in (0,1)$,
$$
\frac{1}{\sqrt{1-\theta}}\leq 2c.
$$
This is obviously impossible as $\theta\To 1$. Thus, we have proved the desired assertion.\vspace{0.2cm}
\end{itemize}
\end{rmk}

The following comparison theorem on unbounded solutions of sub-quadratic BSDEs is one of the main results of the paper. Its proof will be provided in next section.

\begin{thm}\label{thm:3.5}
Assume that $\xi$ and $\xi'$ are two terminal conditions satisfying $\ps, \xi\leq\xi'$, both $g$ and $g'$ are continuous generators in $(y,z)$, $g$ (resp. $g'$) satisfies assumption \ref{A:A3} and
\begin{equation}\label{eq:3.7}
\mathbb{E}\left[\exp\left(p\nu_T^{1-\frac{2}{\delta^\ast}}\left(1+\int_0^T
e^{\mu_s}f_s\mathrm{d}s\right)^{\frac{2}{\delta^\ast}}\right)\right]<+\infty,\ \ \RE p>1.
\end{equation}
Let $(Y_t,Z_t)_{t\in[0,T]}$ and $(Y'_t,Z'_t)_{t\in[0,T]}$ be adapted solutions of BSDE $(\xi,T,g)$ and BSDE $(\xi',T,g')$, respectively, such that
\begin{equation}\label{eq:3.8}
\mathbb{E}\left[\exp\left(p \nu_T^{1-\frac{2}{\delta^\ast}}\sup_{t\in\T}\left( e^{\mu_t}(|Y_t|+|Y'_t|)\right)^{\frac{2}{\delta^\ast}}\right)\right]<+\infty,\ \ \RE p>1.
\end{equation}
If
\begin{equation}\label{eq:3.9}
\as,~~g(t,Y'_t,Z'_t)\leq g'(t,Y'_t,Z'_t)~~~(resp.~~g(t,Y_t,Z_t)\leq g'(t,Y_t,Z_t)),\vspace{0.1cm}
\end{equation}
then $\ps$, for each $t\in[0,T]$, $Y_t\leq Y'_t$.\vspace{0.2cm}

Moreover, if the process $\kappa_\cdot$ in \ref{A:A3} vanishes, then the above assertion remains true when \eqref{eq:3.8} is replaced with the following weaker condition
\begin{equation}\label{eq:3.9*}
\mathbb{E}\left[\sup_{t\in\T}\exp\left(p\nu_t^{1-\frac{2}{\delta^\ast}} \left(e^{\mu_t}(|Y_t|+|Y'_t|)\right)^{\frac{2}{\delta^\ast}}\right)\right]<+\infty,\ \ \RE p>1.
\end{equation}
\end{thm}\vspace{0.01cm}

\begin{rmk}\label{rmk:3.6}
\cref{thm:3.5} strengthens \citet[Theorm 3.8]{FH1}, where $\beta_\cdot$, $\gamma_\cdot$ and $\kappa_\cdot$ are all nonnegative constants and the terminal time $T<+\infty$, and \citet[Theorm 4.2]{GWF} to some extent, where $\beta_\cdot$, $\gamma_\cdot$ and $\kappa_\cdot$ are all nonnegative deterministic functions.
\end{rmk}

\subsection{Existence and uniqueness\vspace{0.2cm}}

According to \cref{thm:2.4} and \cref{thm:3.5}, we can easily derive the following general existence and uniqueness theorem for unbounded solutions of sub-quadratic BSDEs with unbounded stochastic coefficients, which is also one of main results of the paper.

\begin{thm}\label{thm:3.7}
Assume that $\xi$ is a terminal condition, $g$ is a generator which is continuous in $(y,z)$ and satisfies assumptions \ref{A:A1} and \ref{A:A2}, and
\begin{equation}\label{eq:3.15}
\mathbb{E}\left[\exp\left(p\nu_T^{1-\frac{2}{\delta^\ast}}\left(e^{\mu_T}|\xi|+
(\delta^\ast)^{\frac{\delta^\ast}{2}}+\displaystyle\int_0^Te^{\mu_s}f_s\mathrm{d}s
\right)^{\frac{2}{\delta^\ast}}\right)\right]<+\infty,\ \ \RE p>1
\end{equation}
then BSDE $(\xi,T,g)$ admits an adapted solution $(Y_t,Z_t)_{t\in[0,T]}$ such that for each $p>1$,
\begin{align}\label{eq:3.16}
\begin{split}
&\E\left[\sup_{t\in [0,T]}\exp\left(p \nu_t^{1-\frac{2}{\delta^\ast}} (e^{\mu_t}|Y_t|)^{\frac{2}{\delta^\ast}}\right)\right]\\
&\ \ \leq\left(\displaystyle\frac{p}{p-1}\right)^p\mathbb{E}\left[\exp\left(p\nu_T^{1-\frac{2}
{\delta^\ast}}\left(e^{\mu_T}|\xi|+(\delta^\ast)^{\frac{\delta^\ast}{2}}+\displaystyle
\int_0^Te^{\mu_s}f_s\mathrm{d}s\right)^{\frac{2}{\delta^\ast}}\right)\right].
\end{split}
\end{align}

Moreover, we assume that $g$ also satisfies assumption \ref{A:A3}. If either $\kappa_\cdot=0$ or $\nu_T\leq M$ for some constant $M>0$, then the adapted solution satisfying \eqref{eq:3.16} is unique.\vspace{0.2cm}
\end{thm}

\begin{proof}
According to \cref{thm:2.4}, we can conclude that BSDE $(\xi,T,g)$ admits an adapted solution $(Y_t,Z_t)_{t\in[0,T]}$ such that \eqref{eq:2.6*} holds for each $t\in\T$. Then, by Doob's martingale maximum inequality and \eqref{eq:3.15}, we deduce that for each $p>1$,
\begin{align*}
\begin{split}
&\E\left[\sup_{t\in [0,T]}\exp\left(p \nu_t^{1-\frac{2}{\delta^\ast}} (e^{\mu_t}|Y_t|)^{\frac{2}{\delta^\ast}}\right)\right]\\
&\ \ \leq\mathbb{E}\left[\underset{t\in[0,T]}\sup\left\{\mathbb{E}\left[\exp\bigg(
\nu_T^{1-\frac{2}{\delta^\ast}}\bigg(e^{\mu_T}|\xi|+(\delta^\ast)^{\frac{\delta^\ast}{2}}
+\int_0^Te^{\mu_r}f_r\mathrm{d}r\bigg)^{\frac{2}{\delta^\ast}}\bigg)\bigg|\mathcal{F}_t\right]
\right\}^p\right]\\
&\ \ \leq\left(\displaystyle\frac{p}{p-1}\right)^p\mathbb{E}\left[\exp\left(p\nu_T^{1-\frac{2}{\delta^\ast}}
\left(e^{\mu_T}|\xi|+(\delta^\ast)^{\frac{\delta^\ast}{2}}+\displaystyle\int_0^Te^{\mu_s}
f_s\mathrm{d}s\right)^{\frac{2}{\delta^\ast}}\right)\right]<+\infty.
\end{split}
\end{align*}
This is just \eqref{eq:3.16}. Moreover, if \ref{A:A3} holds and either $\kappa_\cdot=0$ or $\nu_T\leq M$ for some $M>0$, then the uniqueness of the adapted solution satisfying \eqref{eq:3.16} follows immediately from \cref{thm:3.5}. The proof is then complete.\vspace{0.2cm}
\end{proof}

\begin{rmk}\label{rmk:3.8}
\cref{thm:3.7} strengthens \citet[Theorm 3.9]{FH1}, where $\beta_\cdot$, $\gamma_\cdot$ and $\kappa_\cdot$ are all nonnegative constants and the terminal time $T<+\infty$, and \citet[Theorm 4.3]{GWF} to some extent, where $\beta_\cdot$, $\gamma_\cdot$ and $\kappa_\cdot$ are all nonnegative deterministic functions.\vspace{0.1cm}
\end{rmk}

\begin{ex}\label{ex:3.9}
The assertions stated in the following several examples can be obtained by our \cref{thm:3.7}, but not by any existing results.
\begin{itemize}
\item [(i)] For each $(\omega,t,y,z)\in \Omega\times\T\times \R\times\R^{d}$, let
$$
g(\omega,t,y,z):=e^{-t}y^2{\bf 1}_{\{y\leq 0\}}+|B_t(\omega)|^{1\over 4}{\bf 1}_{\{0\leq t\leq 1\}}|z|^{3\over 2}+e^{-t}\sin(2\pi |z|){\bf 1}_{\{|z|\leq 1\}}.
$$
According to (i) and (iii) of \cref{pro:3.2} and (v) of \cref{pro:3.3}, it is not very difficult to verify that this $g$ satisfies \ref{A:A1}, \ref{A:A2} and \ref{A:A3} with $f_t:=5e^{-t}$,
$$
\beta_t:=0,\ \ \gamma_t:=|B_t|^{1\over 4}{\bf 1}_{\{0\leq t\leq 1\}}+e^{-t},\ \ \kappa_t:=0,\ \ \delta:=\frac{3}{2}\ \ {\rm and}\ \ \varphi_t(x):=e^{-t}(x+1).
$$
Note that
$$\nu_T=\int_0^T \gamma^4_t\, {\rm d}t\leq 8\int_0^1|B_t|\, {\rm d}t+2,$$
and then $\E[\exp(k\nu_T^{1\over 3})]<+\infty$ for each $k>0$. It then follows from \cref{thm:3.7} that for each $\mathcal{F}_T$-measurable random variable $\xi$ satisfying that
$$\E\left[\exp\left(\left(8\int_0^1|B_t|\, {\rm d}t+2\right)^{1\over 3}|\xi|^{2\over 3}\right)\right]<+\infty,$$
BSDE $(\xi,T,g)$ admits a unique solution $(Y_t,Z_t)_{t\in \T}$ such that for each $p>1$,
$$
\E\left[\sup_{t\in\T}\exp\left(p\nu_t^{1\over 3}|Y_t|^{2\over 3}\right)\right]\leq \E\left[\exp\left(p\left(8\int_0^1|B_t|\, {\rm d}t+2\right)^{1\over 3}\left(|\xi|+\sqrt{27}\,\right)^{2\over 3}\right)\right].
$$

\item [(ii)] Let $T=+\infty$ and $M>0$ be a constant. For each $(\omega,t,y,z)\in \Omega\times\T\times \R\times\R^{d}$, let
$$
g(\omega,t,y,z):=e^{-y}-\sqrt{y}{\bf 1}_{\{y\geq 0\}}+|B_t|^2{\bf 1}_{\{0\leq t\leq \sigma_1\}}|y|+|B_t|{\bf 1}_{\{0\leq t\leq \sigma_2\}}|z|^{6\over 5}+|z|,
$$
where $\sigma_1$ and $\sigma_2$ are two $(\F_t)$-stopping times defined by
$$
\sigma_1:=\inf\left\{t\geq 0: \int_0^t |B_s|^2\, {\rm d}s\geq M\right\}
$$
and
$$
\sigma_2:=\inf\left\{t\geq 0: \int_0^t |B_s|^{5\over 2}\, {\rm d}s\geq M\right\}.\vspace{0.2cm}
$$
According to (i) of \cref{pro:3.2} and (iii) and (v) of \cref{pro:3.3}, it is not very difficult to verify that this $g$ satisfies \ref{A:A1}, \ref{A:A2} and \ref{A:A3} with $f_t:=1$,
$$
\beta_t:=|B_t|^2{\bf 1}_{\{0\leq t\leq \sigma_1\}},\ \gamma_t:=|B_t|{\bf 1}_{\{0\leq t\leq \sigma_2\}}+1,\ \kappa_t:=0,\ \delta:=\frac{6}{5}\ {\rm and}\ \varphi_t(x):=e^x+1+\gamma_t x.
$$
Note by the definitions of $\sigma_1$ and $\sigma_2$ that $\int_0^T (|\beta_t|+\gamma^{5\over 2}_t)\, {\rm d}t\leq 5M+4T$. It then follows from \cref{thm:3.7} that for each $\mathcal{F}_T$-measurable random variable $\xi$ satisfying that
$$\E\left[\exp\left(k|\xi|^{1\over 3}\right)\right]<+\infty,\ \ \RE k>0,$$
BSDE $(\xi,T,g)$ admits a unique solution $(Y_t,Z_t)_{t\in \T}$ such that for each $p>1$,
$$
\E\left[\sup_{t\in\T}\exp\left(pY_t|^{1\over 3}\right)\right]\leq K\E\left[\exp\left(pK|\xi|^{1\over 3}\right)\right],
$$
where $K>0$ is a constant depending only on $M$.

\item [(iii)] Let $T$ be a finite positive real number, and for each $(\omega,t,y,z)\in \Omega\times\T\times \R\times\R^{d}$, let
$$
g(\omega,t,y,z):=\left(\sqrt{-y}-y^3 e^{|B_t|^2}\right){\bf 1}_{\{y<0\}}-\sqrt[3]{y}{\bf 1}_{\{y\geq 0\}}+|B_t| y+e^{\frac{1}{3}\int_0^t |B_s|{\rm d}s}|z|^{4\over 3}.
$$
According to (i) of \cref{pro:3.2} and (ii), (iii) and (v) of \cref{pro:3.3}, it is not very difficult to verify that this $g$ satisfies \ref{A:A1}, \ref{A:A2} and \ref{A:A3} with $f_t:=0$,
$$
\beta_t:=|B_t|,\ \ \gamma_t:=e^{\frac{1}{3}\int_0^t |B_s|{\rm d}s},\ \ \kappa_t:=1,\ \ \delta:=\frac{4}{3}\ \ {\rm and}\ \ \varphi_t(x):=e^{|B_t|^2} x^3+(|B_t|+1)x.
$$
Note that $e^{-(\delta-1)\mu_t}\gamma_t\leq 1$ for each $t\in\T$ and then $\nu_T\leq T$. It then follows from \cref{thm:3.7} that for each $\mathcal{F}_T$-measurable random variable $\xi$ satisfying that
$$\E\left[\exp\left(p\left(e^{\int_0^T|B_t|\,{\rm d}t}|\xi|\right)^{1\over 2}\right)\right]<+\infty,\ \ \RE p>0,$$
BSDE $(\xi,T,g)$ admits a unique solution $(Y_t,Z_t)_{t\in \T}$ such that for each $p>1$,
$$
\E\left[\sup_{t\in\T}\exp\left(p\left(e^{\int_0^t|B_s|\,{\rm d}s}|Y_t|\right)^{1\over 2}\right)\right]\leq K\E\left[\exp\left(pK\left(e^{\int_0^T|B_t|\,{\rm d}t}|\xi|\right)^{1\over 2}\right)\right],
$$
where $K>0$ is a constant depending only on $T$.

\item [(iv)] Let $d=1$, $T<+\infty$ and for each $(\omega,t,y,z)\in \Omega\times\T\times \R\times\R^{d}$,
$$
g(\omega,t,y,z):=B_t^5 y +y^4 {\bf 1}_{\{y\leq 0\}}+e^{\frac{1}{4}\int_0^t B_s^5{\rm d}s}|z|^{5\over 4}.
$$
According to (i) of \cref{pro:3.2} and (ii) and (v) of \cref{pro:3.3}, it is not very difficult to verify that this $g$ satisfies \ref{A:A1}, \ref{A:A2} and \ref{A:A3} with $f_t:=0$,
$$
\beta_t:=B_t^5,\ \ \gamma_t:=e^{\frac{1}{4}\int_0^t B_s^5{\rm d}s},\ \ \kappa_t:=0,\ \ \delta:=\frac{5}{4}\ \ {\rm and}\ \ \varphi_t(x):=|B_t|^5 x+x^4.
$$
Note that $e^{-(\delta-1)\mu_t}\gamma_t\leq 1$ for each $t\in\T$ and then $\nu_T\leq T$. It then follows from \cref{thm:3.7} that for random variable $\xi:=e^{-\int_0^T B_t^5\, {\rm d}t}$, BSDE $(\xi,T,g)$ admits a unique adapted solution $(Y_t,Z_t)_{t\in \T}$ such that
$$
|Y_t|\leq K e^{-\int_0^t B_s^5\, {\rm d}s},\ \ t\in\T,\vspace{-0.2cm}
$$
where $K>0$ is a constant depending only on $T$. We would like to especially mention that $\beta_\cdot$ can take negative values and that this $\xi$ is not an integrable random variable.\vspace{0.2cm}
\end{itemize}
\end{ex}

\section{Proofs of some theorems and propositions}
\label{sec:3-Stability theorem}
\setcounter{equation}{0}

\subsection{The proof of \cref{pro:2.2}\vspace{0.2cm}}

In this subsection, we give the proof of \cref{pro:2.2}.

\begin{proof}[Proof of \cref{pro:2.2}]
We first consider the case of $\beta_\cdot\equiv 0$, and then the general case.\vspace{0.2cm}

{\bf Case 1}. The case of $\beta_\cdot\equiv 0$. We first define a function
$$\psi(x,y):=\exp(xy^{\frac{2}{\delta^\ast}}),\ \ (x,y)\in[1,+\infty)\times[y_0,+\infty),$$
where $y_0$ is a constant bigger than $(\delta^\ast)^{\frac{\delta^\ast}{2}}$. Note that the larger the $\delta^\ast$, the smaller the $\psi(x,y)$. For each $(x,y)\in[1,+\infty)\times[y_0,+\infty)$, a simple calculation yields that
\begin{align}\label{eq:2.7}
\begin{split}
\psi_x(x,y)&=\psi(x,y)y^{\frac{2}{\delta^\ast}}>0,\\
\psi_y(x,y)&=\psi(x,y)x\frac{2}{\delta^\ast}y^{{\frac{2}{\delta^\ast}}-1}
=\psi_x(x,y)\frac{2}{\delta^\ast}\frac{x}{y}>0,\\
\psi_{yy}(x,y)&=\psi(x,y)x\frac{2}{\delta^\ast}y^{{\frac{2}{\delta^\ast}}-2}
\left[\frac{2}{\delta^\ast}xy^{\frac{2}{\delta^\ast}}-(1-\frac{2}{\delta^\ast})\right]\\
&\geq \psi(x,y)x\frac{2}{\delta^\ast}y^{{\frac{2}{\delta^\ast}}-2}\left[\frac{2}
{\delta^\ast}y_0^{\frac{2}{\delta^\ast}}-1\right]>0,
\end{split}
\end{align}
and in light of $\frac{1}{\delta^\ast}xy^{\frac{2}{\delta^\ast}}\geq\frac{1}{\delta^\ast}
[(\delta^\ast)^{\frac{\delta^\ast}{2}}]^{\frac{2}{\delta^\ast}}=1$, we have
\begin{align}\label{eq:2.8}
\psi_{yy}(x,y)&\geq\psi_x(x,y)\frac{2}{(\delta^\ast)^2}x^2y^{\frac{2}{\delta^\ast}-2}.
\end{align}
In the sequel, we define a process
\begin{align}\label{eq:2.8-1}
\bar{Y}_s:=|Y_s|+(\delta^\ast)^{\frac{\delta^\ast}{2}}+\int_0^sf_r\mathrm{d}r,\ \ s\in \T.
\end{align}
It follows from It\^{o}-Tanaka's formula and assumption \ref{A:A1} with $\beta_\cdot\equiv 0$ that
\begin{align}\label{eq:2.9}
\begin{split}
 \mathrm{d}\bar{Y}_s&=\left[-\mathrm{sgn}(Y_s)g(s,Y_s,Z_s)+f_s\right]\mathrm{d}s+\mathrm{sgn}(Y_s)Z_s\cdot
 \mathrm{d}B_s+\mathrm{d}L_s\\
 &\geq -\gamma_s|Z_s|^\delta\mathrm{d}s+\mathrm{sgn}(Y_s)Z_s\cdot
 \mathrm{d}B_s,\ \ s\in[0,T],
\end{split}
\end{align}
where $L(\cdot)$ represent the local time of $Y_\cdot$ at time 0. Then, we apply It\^{o} formula to $\psi(\bar\nu_s^{1-\frac{2}{\delta^\ast}},\bar{Y}_s)$, and combine with \eqref{eq:2.7}, \eqref{eq:2.8} and \eqref{eq:2.9} to obtain that
\begin{align}\label{eq:2.10}
 \begin{split}
&\mathrm{d}\psi(\bar\nu_s^{1-\frac{2}{\delta^\ast}},\bar{Y}_s)\\
&\ \ =\psi_x(\bar\nu_s^{1-\frac{2}
{\delta^\ast}},\bar{Y}_s)\mathrm{d}\bar\nu_s^{1-\frac{2}{\delta^\ast}}+
\psi_y(\bar\nu_s^{1-\frac{2}{\delta^\ast}},\bar{Y}_s)\mathrm{d}\bar{Y}_s+
\frac{1}{2}\psi_{yy}(\bar\nu_s^{1-\frac{2}{\delta^\ast}},\bar{Y}_s)\mathrm{d}\langle\bar{Y}\rangle_s\\
&\ \ \geq \psi_x(\bar\nu_s^{1-\frac{2}{\delta^\ast}},\bar{Y}_s)\mathrm{d}\bar\nu_s^{1-\frac{2}
 {\delta^\ast}}+\left(-\gamma_s\psi_y(\bar\nu_s^{1-\frac{2}{\delta^\ast}},\bar{Y}_s)|Z_s|^\delta
 +\frac{1}{2}\psi_{yy}(\bar\nu_s^{1-\frac{2}{\delta^\ast}},\bar{Y}_s)|Z_s|^2\right)\mathrm{d}s\\
&\ \ \ \ +\psi_y(\bar\nu_s^{1-\frac{2}{\delta^\ast}},\bar{Y}_s)\mathrm{sgn}(Y_s)Z_s\cdot\mathrm{d}B_s\\
&\ \ \geq \psi_x(\bar\nu_s^{1-\frac{2}{\delta^\ast}},\bar{Y}_s)\left[\mathrm{d}
 \bar\nu_s^{1-\frac{2}{\delta^\ast}}+\frac{1}{\delta^\ast}\frac{
 \bar\nu_s^{1-\frac{2}{\delta^\ast}}}{\bar{Y}_s^2}\left(-2\gamma_s\bar{Y}_s|Z_s|^\delta
 +\frac{\bar\nu_s^{1-\frac{2}{\delta^\ast}}
 }{\delta^\ast}
 \bar{Y}_s^{\frac{2}{\delta^\ast}}|Z_s|^2\right)\right]\mathrm{d}s\\
&\ \ \ \ +\psi_y(\bar\nu_s^{1-\frac{2}{\delta^\ast}},\bar{Y}_s)\mathrm{sgn}(Y_s)Z_s\cdot\mathrm{d}B_s,
\end{split}
\end{align}
where $\langle\bar{Y}\rangle$ represents the quadratic variation process of $\bar{Y}$. In light of \eqref{eq:2.4}, it follows from Young's inequality that for each $s\in[0,T]$,
\begin{align}\label{eq:2.11}
\begin{split}
2\gamma_s\bar{Y}_s|Z_s|^\delta&=\left(\frac{2}{\delta}\cdot\frac{\bar\nu_s^{1-\frac{2}
{\delta^\ast}}}{\delta^\ast}\bar{Y}_s^{\frac{2}{\delta^\ast}}|Z_s|^2\right)^{\frac{\delta}{2}}
\left(2\left(\delta\delta^\ast\right)
^{\frac{\delta}{2-\delta}}\gamma_s^{\frac{2}{2-\delta}}\bar{Y}_s^2\bar\nu_s^{-1}\right)^{1-\frac{\delta}{2}}\\
&\leq\frac{\bar\nu_s^{1-\frac{2}{\delta^\ast}}}{\delta^\ast}\bar{Y}^{\frac{2}{\delta^\ast}}_s|Z_s|^2+
\frac{2-\delta}{\delta-1}C_{\delta}\gamma_s^{\frac{2}{2-\delta}}\bar{Y}_s^2\bar\nu_s^{-1},
\end{split}
\end{align}
where the constant $C_\delta>0$ is defined in \eqref{eq:2.3}. Combining with \eqref{eq:2.4} and \eqref{eq:2.11} yields that
\begin{align}\label{eq:2.12}
\begin{split}
&\mathrm{d}\bar\nu_s^{1-\frac{2}{\delta^\ast}}+\frac{1}{\delta^\ast}
\frac{\bar\nu_s^{1-\frac{2}{\delta^\ast}}}{\bar{Y}_s^2}\left(-2\gamma_s
\bar{Y}_s|Z_s|^\delta+\frac{\bar\nu_s^{1-\frac{2}{\delta^\ast}}}{\delta^\ast}
 \bar{Y}_s^{\frac{2}{\delta^\ast}}|Z_s|^2\right)\mathrm{d}s\\
&\ \ \geq (1-\frac{2}{\delta^\ast})\bar\nu_s^{-\frac{2}{\delta^\ast}}
 \mathrm{d}\bar\nu_s+\frac{\bar\nu_s^{-\frac{2}{\delta^*}}}{\delta^*}
\frac{2-\delta}{\delta-1}C_\delta\gamma_s^{\frac{2}{2-\delta}}\mathrm{d}s=0.
\end{split}
\end{align}
Furthermore, for each $t\in[0,T]$ and integer $m\geq1$, let us denote the following stopping time
$$\tau_m^t:=\inf\displaystyle\left\{s\in[t,T]:\int_t^s(\psi_y
(\bar\nu_s^{1-\frac{2}{\delta^\ast}},\bar{Y}_s))^2|Z_r|^2\mathrm{d}r\geq m\right\}\wedge T$$
with the convention $\inf\emptyset=+\infty$. In light of \eqref{eq:2.12}, by integrating on the interval $[t,\tau_m^t]$ on both sides of \eqref{eq:2.10} we can obtain that for each $m\geq1$,
$$
\exp(\bar\nu_t^{1-\frac{2}{\delta^\ast}}\bar{Y}_t^{\frac{2}{\delta^\ast}})
\leq\exp(\bar\nu_{\tau_m^t}^{1-\frac{2}{\delta^\ast}}\bar{Y}_{\tau_m^t}^{\frac{2}{\delta^\ast}})
-\int_t^{\tau_m^t}\psi_y(\bar\nu_s^{1-\frac{2}{\delta^\ast}},\bar{Y}_s)\mathrm{sgn}(Y_s)Z_s\cdot\mathrm{d}B_s.
$$
By taking the conditional expectation in the last inequality, we get that for each $m\geq1$,
\begin{equation}\label{eq:2.13}
\exp(\bar\nu_t^{1-\frac{2}{\delta^\ast}}\bar{Y}_t^{\frac{2}{\delta^\ast}})
\leq\mathbb{E}\left[\exp(\bar\nu_{\tau_m^t}^{1-\frac{2}{\delta^\ast}}
\bar{Y}_{\tau_m^t}^{\frac{2}{\delta^\ast}})\bigg|\mathcal{F}_t\right],\ \ t\in [0,T].
\end{equation}
Finally, in light of Fatou's lemma and Lebesgue's dominated convergence theorem along with \eqref{eq:2.5} with $\mu_\cdot=0$, sending $m$ to infinity in \eqref{eq:2.13} yields that for each $t\in[0,T]$,
\begin{equation}\label{eq:2.14}
\exp(\bar\nu_t^{1-\frac{2}{\delta^\ast}}|Y_t|^{\frac{2}{\delta^\ast}})\leq
\mathbb{E}\left[\exp\left(\bar\nu_T^{1-\frac{2}{\delta^\ast}}\left(|\xi|
+(\delta^\ast)^{\frac{\delta^\ast}{2}}+\int_0^Tf_s\mathrm{d}s
\right)^{\frac{2}{\delta^\ast}}\right)\bigg|\mathcal{F}_t\right].
\end{equation}
In light of \eqref{eq:2.4}, the desired inequality \eqref{eq:2.6} follows immediately from \eqref{eq:2.14} for this case.\vspace{0.2cm}

{\bf Case 2}. The general case of  $\beta_\cdot\not\equiv 0$. For each $t\in[0,T]$, let
$$\tilde{Y}_t:=e^{\mu_t}Y_t,~~\tilde{Z}_t:=e^{\mu_t}Z_t~~\mathrm{and}~~\tilde{\xi}=e^{\mu_T}\xi.$$
Then, in light of \eqref{eq:2.2}, the pair $(\tilde{Y}_t,\tilde{Z}_t)_{t\in\T}$ solves the following BSDE
\begin{equation}\label{eq:2.15}
\tilde{Y}_t=\tilde{\xi}+\int_{t}^{T} \tilde{g}(s,\tilde{Y}_s,\tilde{Z}_s)\mathrm{d}s-\int_{t}^{T}\tilde{Z}_s\cdot \mathrm{d}B_s,\ \ t\in [0,T],
\end{equation}
where $\ass$, for each $(s,y,z)\in [0,T]\times\mathbb{R}\times\mathbb{R}^d$,
\begin{equation}\label{eq:2.16}
\tilde{g}(s,y,z):=e^{\mu_s}g(s,e^{-\mu_s}y,e^{-\mu_s}z)-\beta_sy.
\end{equation}
It follows from \ref{A:A1} that $\ass$, for each $(y,z)\in \mathbb{R}\times\mathbb{R}^d$,
\begin{equation}\label{eq:2.17}
\mathrm{sgn}(y)\tilde{g}(\omega,s,y,z)\leq\tilde{f}_s(\omega)+\tilde{\gamma}_s(\omega)|z|^\delta,
\end{equation}
where
$$\tilde{f}_s(\omega):=e^{\mu_s}f_s(\omega)\ \ {\rm and}\ \  \tilde{\gamma}_s(\omega):=e^{-(\delta-1)\mu_s}\gamma_s(\omega).
$$
Consequently, the generator $\tilde{g}$ satisfies \ref{A:A1} with $f_s:=\tilde{f}_s$, $\beta_s:=0$ and $\gamma_s:=\tilde{\gamma}_s$. Thus, by applying \eqref{eq:2.14} obtained in Step 1 to BSDE \eqref{eq:2.15}, the desired assertion \eqref{eq:2.6} follows immediately provided that \eqref{eq:2.5} is in force. The proof of \cref{pro:2.2} is complete.\vspace{0.2cm}
\end{proof}

\subsection{The proof of \cref{thm:2.4}\vspace{0.2cm}}

In this subsection, we give the proof of \cref{thm:2.4}.

\begin{proof}[Proof of \cref{thm:2.4}]
We first consider the case that $\beta_\cdot\geq 0$ and $g$ satisfies \ref{A:A1} and \ref{A:A2'}, and then consider the other case that $\beta_\cdot\in \R$ and $g$ satisfies \ref{A:A1} and \ref{A:A2}.\vspace{0.2cm}

{\bf Case 1}. The case that $\beta_\cdot\geq 0$ and $g$ satisfies \ref{A:A1} and \ref{A:A2'}. For each pair of positive integers $n,q\geq 1$ and each $(\omega,t,y,z)\in\Omega\times[0,T]\times\mathbb{R}\times\mathbb{R}^d$, let
$$\xi^{n,q}:=\xi^+\wedge (ne^{-\mu_T})-\xi^-\wedge (qe^{-\mu_T})$$
and
$$g^{n,q}(\omega,t,y,z):=g^+(\omega,t,y,z)\wedge (ne^{-t-\mu_t})-g^-(\omega,t,y,z)\wedge (qe^{-t-\mu_t}).$$
It is evident that $|\xi^{n,q}|\leq|\xi|\wedge [(n\vee q)e^{-\mu_T}]$ and $|g^{n,q}|\leq |g|\wedge[(n\vee q)e^{-t-\mu_t}]$. It is also easy to verify that $\as$, $g^{n,q}$ is continuous in $(y,z)$ and that $g^{n,q}$ satisfies assumptions \ref{A:A2'} and \ref{A:A1} with $f_t$ being replaced with $f_t^{n,q}:=f_t\wedge[(n\vee q)e^{-t-\mu_t}]$ satisfying
\begin{equation}\label{eq:2.20}
\int_0^T e^{\mu_s}f_s^{n,q}\mathrm{d}s\leq(n\vee q)\int_0^Te^{-s}\mathrm{d}s\leq n\vee q.
\end{equation}
Note that $|g^{n,q}(t,y,z)|\leq u_t:=(n\vee q)e^{-t}$ for each $(y,z)\in\mathbb{R}\times\mathbb{R}^d$ with $\int_0^Tu_tdt\leq n\vee q$, and that $\xi^{n,q}\in L^\infty(\Omega,\mathcal{F}_T,\mathbb{P})$. By \citet[Lemma 4]{FL} we conclude that the following BSDE admits a minimal bounded solution $(Y_t^{n,q},Z_t^{n,q})_{t\in[0,T]}\in \mathcal{S}^\infty\times\mathcal{M}^2$:
\begin{equation}\label{eq:2.22}
Y_t^{n,q}=\xi^{n,q}+\int_{t}^{T} g^{n,q}(s,Y_s^{n,q},Z_s^{n,q})\mathrm{d}s-\int_{t}^{T}Z_s^{n,q}\cdot \mathrm{d}B_s,~~~~~t\in[0,T].
\end{equation}
By virtue of the comparison result of \citet[Theorem 2]{FL}, $Y_\cdot^{n,q}$ is nondecreasing in $n$ and non-increasing in $q$. Furthermore, by virtue of the definitions of $\xi^{n,q}$ and $g^{n,q}$, by taking the conditional mathematical expectation with respect to $\mathcal{F}_t$ in \eqref{eq:2.22} we can obtain
\begin{align*}
\begin{split}
e^{\mu_t}|Y_t^{n,q}|&\leq\mathbb{E}\left[e^{\mu_t}|\xi^{n,q}|+\int_t^T e^{\mu_t}|g^{n,q}(s,Y_s^{n,q},Z_s^{n,q})|
\mathrm{d}s\bigg|\mathcal{F}_t\right]\\
&\leq \mathbb{E}\left[(n\vee q)e^{-\mu_T+\mu_t}+(n\vee q)\int_t^Te^{-s-\mu_s+\mu_t}\mathrm{d}s\bigg|\mathcal{F}_t\right]\leq 2(n\vee q),~~~t\in[0,T],
\end{split}
\end{align*}
which along with \eqref{eq:2.18} and \eqref{eq:2.20} implies that the condition \eqref{eq:2.5} is fulfilled for $(Y^{n,q}_\cdot,f^{n,q}_\cdot)$. It then follows from \cref{pro:2.2} that for each $n,q\geq 1$ and $t\in [0,T]$,
\begin{align}\label{eq:2.23}
\begin{split}
\exp\left(|Y_t^{n,q}|^{\frac{2}{\delta^\ast}}\right)&\leq \exp\left(\nu_t^{1-\frac{2}{\delta^\ast}}(e^{\mu_t}|Y_t^{n,q}|)^{\frac{2}{\delta^\ast}}\right)\\
&\leq \mathbb{E}\left[\exp\left(\nu_T^{1-\frac{2}{\delta^\ast}}\left(e^{\mu_T}|\xi^{n,q}|
+(\delta^\ast)^{\frac{\delta^\ast}{2}}+\displaystyle\int_0^T e^{\mu_s}f_s^{n,q}\mathrm{d}s
\right)^{\frac{2}{\delta^\ast}}\right)\bigg|\mathcal{F}_t\right]\\
&\leq\mathbb{E}\left[\exp\left(\nu_T^{1-\frac{2}{\delta^\ast}}\left(e^{\mu_T}|\xi|
+(\delta^\ast)^{\frac{\delta^\ast}{2}}+\displaystyle\int_0^T e^{\mu_s}f_s\mathrm{d}s
\right)^{\frac{2}{\delta^\ast}}\right)
\bigg|\mathcal{F}_t\right]=:H_t.
\end{split}
\end{align}
In light of \eqref{eq:2.19}, the process $(H_t)_{t\in\T}$ is well defined.\vspace{0.2cm}

Now, we can apply the localization procedure developed initially in \citet{BH}. For each pair of integers $m,l\geq1$, we define the following stopping times:
$$\sigma_m:=\inf\left\{t\in[0,T]:H_t\geq\exp(m^{\frac{2}{\delta^\ast}})\right\}\wedge T$$
and
$$
\tilde{\sigma}_{m,l}:=\inf\left\{t\in[0,T]:\int_0^t \bar\varphi_s(m){\rm d}s\geq l\right\}\wedge \sigma_m\vspace{0.2cm}
$$
with the convention $\inf\emptyset=+\infty$, where the function $\bar\varphi_t(x)$
is defined in assumption \ref{A:A2'}. Then, the pair $(Y_{m,l}^{n,q}(t),Z_{m,l}^{n,q}(t)):=(Y^{n,q}_{t\wedge \tilde{\sigma}_{m,l}},Z_t^{n,q}\mathbf{1}_{t\leq\tilde{\sigma}_{m,l}})$ is a solution to the following BSDE:
\begin{equation*}
Y_{m,l}^{n,q}(t)=Y_{\tilde{\sigma}_{m,l}}^{n,q}+\int_{t}^{T} \mathbf{1}_{s\leq\tilde{\sigma}_{m,l}}g^{n,q}(s,Y_{m,l}^{n,q}(s),Z_{m,l}^{n,q}(s))
\mathrm{d}s-\int_{t}^{T}Z_{m,l}^{n,q}(s)\cdot \mathrm{d}B_s,~~~~~t\in[0,T].
\end{equation*}
Note that for any fixed $m,l\geq1$, $Y_{m,l}^{n,q}(\cdot)$ is nondecreasing in $n$ and non-increasing in $q$, and that the sequence $(g^{n,q})_{n,q}$ converges locally uniformly in $(y,z)$ to the generator $g$ as $n,q\rightarrow\infty$. In light of \eqref{eq:2.23} along with the definitions of $\sigma_m$ and $\tilde{\sigma}_{m,l}$, we have
$$\as,\ \ \underset{n,q\geq1}\sup|Y_{m,l}^{n,q}(t)|\leq m.$$
Furthermore, since $|g^{n,q}|\leq |g|$, it follows from \ref{A:A2'} and the definition of $\tilde{\sigma}_{m,l}$ that $\ass$,
\begin{align*}
\begin{split}
\forall(y,z)\in[-m,m]\times\mathbb{R}^d,~\underset{n,q\geq1}\sup(
\mathbf{1}_{s\leq\tilde{\sigma}_{m,l}}|g^{n,q}(s,y,z)|)&\leq \mathbf{1}_{s\leq\tilde{\sigma}_{m,l}}\bar\varphi_s(m)+\phi(m)|z|^2,
\end{split}
\end{align*}
with $\int_0^T\mathbf{1}_{s\leq\tilde\sigma_{m,l}}\bar\varphi_s(m)\leq l$. Thus, we can apply the stability property on the bounded solutions of BSDEs (see for example \citet[Proposition 3.1]{LF}) to conclude that the pair $(Y_{m,l}(t),Z_{m,l}(t))_{t\in\T}:=(\underset{q}\inf ~\underset{n}\sup Y^{n,q}_{t\wedge \tilde{\sigma}_{m,l}},\underset{n,q\rightarrow\infty}\lim Z_t^{n,q}\mathbf{1}_{t\leq\tilde{\sigma}_{m,l}})_{t\in\T}$ solves the following BSDE:
$$Y_{m,l}(t)=\underset{q}\inf ~\underset{n}\sup Y^{n,q}_{\tilde{\sigma}_{m,l}}+\int_t^Tg(s,Y_{m,l}(s),Z_{m,l}(s))\mathrm{d}s-
\int_t^TZ_{m,l}(s)\cdot\mathrm{d}B_s,~~t\in[0,T],$$
where $Y_{m,l}(\cdot)$ is continuous and the limit of sequence $(Z_t^{n,q}\mathbf{1}_{t\leq\tilde{\sigma}_{m,l}})_{t\in [0,T]}$ is taken in $\mathcal{M}^2$. Finally, in light of the last equation along with the stability of stopping times $\sigma_m$ and $\tilde{\sigma}_{m,l}$, combining with the facts that $\mathrm{d}\mathbb{P}\times\mathrm{d}t-a.e.,$ for any $m,l\geq 1$,
$$Z_{m+1,l+1}\mathbf{1}_{t\leq\tilde{\sigma}_{m,l}}=Z_{m,l+1}\mathbf{1}_{t\leq\tilde{\sigma}_{m,l}}
=Z_{m,l}\mathbf{1}_{t\leq\tilde{\sigma}_{m,l}}
=\underset{n,q\rightarrow\infty}\lim Z_t^{n,q}\mathbf{1}_{t\leq\tilde{\sigma}_{m,l}}$$
and
$$Y_{m+1,l+1}(t\wedge\tilde{\sigma}_{m,l})=Y_{m,l+1}(t\wedge\tilde{\sigma}_{m,l})
=Y_{m,l}(t\wedge\tilde{\sigma}_{m,l})=\underset{q}\inf ~\underset{n}\sup Y^{n,q}_{t\wedge \tilde{\sigma}_{m,l}},$$
we can obtain an adapted solution $(Y_t,Z_t)_{t\in[0,T]}$ of BSDE $(\xi,T,g)$, where
$$
Y_t:=\underset{q}\inf ~\underset{n}\sup Y^{n,q}_{t} ~\mathrm{and} ~ Z_t:=\sum_{m=1}^{+\infty}\left(\sum_{l=1}^{+\infty}Z_{m,l}(t)\mathbf{1}_{t\in
[\tilde{\sigma}_{m,l-1},\tilde{\sigma}_{m,l})}\right)
\mathbf{1}_{t\in[\sigma_{m-1},\sigma_m)},~~t\in[0,T],
$$
where $\sigma_0:=0$ and $\tilde{\sigma}_{m,0}:=0$ for each $m\geq1$. And, for the case of $\beta_\cdot\geq 0$, the desired assertion \eqref{eq:2.6*} follows immediately from \eqref{eq:2.23} by sending $n,q\rightarrow\infty$.\vspace{0.2cm}

{\bf Case 2}. The case that $\beta_\cdot\in \R$ and $g$ satisfies assumptions \ref{A:A1} and \ref{A:A2}. First, we consider the following BSDE
\begin{equation}\label{eq:2.15*}
\tilde{Y}_t=\tilde{\xi}+\int_{t}^{T} \tilde{g}(s,\tilde{Y}_s,\tilde{Z}_s)\mathrm{d}s-\int_{t}^{T}\tilde{Z}_s\cdot \mathrm{d}B_s,\ \ t\in [0,T],
\end{equation}
where $\tilde{\xi}=e^{\mu_T}\xi$ and
\begin{equation}\label{eq:2.24*}
\tilde{g}(\omega,t,y,z):=e^{\mu_t}g(\omega,t,e^{-\mu_t}y,e^{-\mu_t}z)-\beta_ty,\ \ (\omega,t,y,z)\in\Omega\times\T\times\R\times\R^d.
\end{equation}
By an identical argument to \eqref{eq:2.17}, we deduce that the generator $\tilde{g}$ satisfies assumption \ref{A:A1} with $(f_t,\beta_t,\gamma_t)$ being replaced with $(\tilde f_t,\tilde \beta_t,\tilde\gamma_t):=(e^{\mu_t}f_t,0,e^{-(\delta-1)\mu_t}\gamma_t)$. And, by assumption \ref{A:A2} of $g$ and Young's inequality we have that $\as$, for each $(y,z)\in\R\times\R^d$,
\begin{align}\label{eq:2.24}
\begin{split}
|\tilde{g}(t,y,z)|&\leq e^{\mu_t}|g(t,e^{-\mu_t}y,e^{-\mu_t}z)|+|\beta_t||y|\\
&\leq e^{\mu_t}\varphi_t(e^{-\mu_t}|y|)+|\beta_t||y|
+\gamma_te^{-(\delta-1)\mu_t}|z|^{\delta}\\
&\leq \tilde\varphi_t(|y|)+\gamma_te^{-(\delta-1)\mu_t}|z|^{\delta},
\end{split}
\end{align}
where for each $(\omega,t,x)\in\Omega\times\T\times\R_+$,
$$
\tilde\varphi_t(x):=e^{\mu_t}\varphi_t(e^{-\mu_t}x)+|\beta_t|x\in \bs.\vspace{-0.1cm}
$$
Consequently, the generator $\tilde{g}$ also satisfies assumption \ref{A:A2} with $(\varphi_t(x),\gamma_t)$ being replaced with $(\tilde\varphi_t(x),e^{-(\delta-1)\mu_t}\gamma_t)$. Furthermore, by an identical argument to \cref{rmk:2.1}, we know that the generator $\tilde g$ also satisfies assumption \ref{A:A2'}. Then, according to Case 1, we can deduce that BSDE \eqref{eq:2.15*} admits a solution $(\tilde Y_t,\tilde Z_t)_{t\in[0,T]}$ such that for each $t\in [0,T]$,
$$
\exp\left(\nu_t^{1-\frac{2}{\delta^\ast}}|\tilde Y_t|^{\frac{2}{\delta^\ast}}\right)\leq \mathbb{E}
\left[\exp\left(\nu_T^{1-\frac{2}{\delta^\ast}}\left(|\tilde\xi|
+(\delta^\ast)^{\frac{\delta^\ast}{2}}+\int_0^T \tilde f_s\mathrm{d}s
\right)^{\frac{2}{\delta^\ast}}\right)
\bigg|\mathcal{F}_t\right].
$$
With the above assertion in hands, by a simple calculation we can conclude that the pair of processes $(Y_t,Z_t)_{t\in\T}:=(e^{-\mu_t}\tilde{Y_t},e^{-\mu_t}\tilde{Z_t})_{t\in\T}$ is an adapted solution of BSDE $(\xi,T,g)$ such that $Y_\cdot$ satisfies \eqref{eq:2.6*}. The proof of \cref{thm:2.4} is then complete. \vspace{0.2cm}
\end{proof}

\begin{rmk}\label{rmk:2.5*}
The method of proof used to address Case 1 cannot be applied to the proof for Case 2. Indeed, for Case 2, since $\beta_\cdot\in\R$, the truncated function $g^{n,q}$ defined in Case 1 can not be proved to satisfy \ref{A:A1} in general, and hence the uniform estimate of $y_\cdot^{n,q}$ in \eqref{eq:2.23} cannot be obtained via \cref{pro:2.2}. Instead, \ref{A:A2'} is strengthened \ref{A:A2} and an exponential shift transformation is applied to transform BSDE $(\xi,T,g)$ into BSDE \eqref{eq:2.15*} whose generator $\tilde g$ satisfies \ref{A:A2'} and \ref{A:A1} with $\beta_\cdot=0$ such that the conclusion of Case 1 can be available. In particular, \ref{A:A2'} is not enough for the case of $\beta_\cdot\in\R$. \vspace{0.2cm}
\end{rmk}

\subsection{The proof of \cref{pro:3.3}\vspace{0.2cm}}

In this subsection, we give the proof of \cref{pro:3.3}.

\begin{proof}[Proof of \cref{pro:3.3}]
(i), (ii), (iv) and (v) are all evident. We only need to prove (iii). \vspace{0.2cm}

We first suppose that $\as$, $g(\omega, t,\cdot)$ is both convex on $\R_-$ and $\R_+$, and then verify that \eqref{eq:3.4} holds. Given arbitrary $(y_1,y_2)\in\R\times\R$ and $\theta\in (0,1)$. Note that \eqref{eq:3.4} is obviously true when $y_1<0$ and $y_2>0$. It suffices to consider the following three cases:\vspace{0.2cm}

(i) The case of $y_1, y_2\geq 0$. In light of the convexity of $g(\omega,t,\cdot)$ on $\R_+$ and the fact of $g(\omega,t,0)=0$, we have $\as$,
$$g(\omega,t,\theta y_2)=g(\omega,t,\theta y_2+(1-\theta)\cdot0)\leq\theta g(\omega,t,y_2).$$
Furthermore, since $g(\omega,t,\cdot)$ is non-increasing, by the last inequality we have $\as$,
$$
\mathbf{1}_{\{y_1>\theta y_2\}}(g(\omega,t,y_1)-\theta g(\omega,t,y_2))\leq \mathbf{1}_{\{y_1>\theta y_2\}}(g(\omega,t,\theta y_2)-\theta g(\omega,t,y_2))\leq 0.
$$

(ii) The case of $y_1, y_2\leq 0$. In light of the convexity and the non-increasing property of $g(\omega,t,\cdot)$ on $\R_-$ and $g(\omega,t,0)=0$, an identical argument to case (i) yields that $\as$,
$$
\mathbf{1}_{\{y_1>\theta y_2\}}(g(\omega,t,y_1)-\theta g(\omega,t,y_2))\leq \mathbf{1}_{\{y_1>\theta y_2\}}(g(\omega,t,\theta y_2)-\theta g(\omega,t,y_2))\leq 0.
$$

(iii) The case of $y_1>0$ and $y_2<0$. In light of the non-increasing property of $g(\omega,t,\cdot)$ on $\R$ and the fact of $g(\omega,t,0)=0$, we have $\as$,
$$
\mathbf{1}_{\{y_1>\theta y_2\}}(f(\omega,t,y_1)-\theta f(\omega,t,y_2))\leq f(\omega,t,0)-\theta f(\omega,t,0)=0.
$$
Consequently, the desired inequality \eqref{eq:3.4} holds for all cases.\vspace{0.2cm}

Next, we suppose that $\as$, $g(\omega, t,y)\leq \kappa_t(\omega)|y|$ for each $y\in\R_-$ and $g(\omega, t,\cdot)$ is convex on $\R_+$, and then verify that \eqref{eq:3.5} holds. Given arbitrary $(y_1,y_2)\in\R\times\R$ and $\theta\in (0,1)$. Note that \eqref{eq:3.5} is true when $y_1<0$ and $y_2>0$. It suffices to consider three cases below:\vspace{0.2cm}

(iv) The case of $y_1, y_2\geq 0$. In light of the convexity and the non-increasing property of $g(\omega,t,\cdot)$ on $\R_+$ and $g(\omega,t,0)=0$, an identical argument to case (i) yields that $\as$,
$$
\mathbf{1}_{\{y_1>\theta y_2\}}(g(\omega,t,y_1)-\theta g(\omega,t,y_2))\leq \mathbf{1}_{\{y_1>\theta y_2\}}(g(\omega,t,\theta y_2)-\theta g(\omega,t,y_2))\leq 0.
$$

(v) The case of $y_1, y_2\leq 0$. In light of the non-increasing property and the upper bound property of $g(\omega,t,\cdot)$ on $\R_-$, we can deduce that $\as$,
\begin{align*}
\begin{split}
\mathbf{1}_{\{y_1>\theta y_2\}}(g(\omega,t,y_1)-\theta g(\omega,t,y_2))
&\leq \mathbf{1}_{\{y_1>\theta y_2\}}(g(\omega,t,y_2)-\theta g(\omega,t,y_2))\\
&=(1-\theta) g(\omega,t,y_2)\leq (1-\theta)\kappa_t(\omega)|y_2|.
\end{split}
\end{align*}

(vi) The case of $y_1>0$ and $y_2<0$. In light of the non-increasing property of $g(\omega,t,\cdot)$ on $\R$ and the fact of $g(\omega,t,0)=0$, we have $\as$,
$$
\mathbf{1}_{\{y_1>\theta y_2\}}(f(\omega,t,y_1)-\theta f(\omega,t,y_2))\leq f(\omega,t,0)-\theta f(\omega,t,0)=0.
$$
Consequently, the desired inequality \eqref{eq:3.5} holds for all cases. The proof is then complete.\vspace{0.2cm}
\end{proof}

\subsection{The proof of \cref{thm:3.5}\vspace{0.2cm}}

In this subsection, we give the proof of \cref{thm:3.5}.

\begin{proof}[Proof of \cref{thm:3.5}]
Let us first assume that \eqref{eq:3.7} and \eqref{eq:3.8} holds, and prove that
\begin{equation}\label{eq:3.10}
\mathbb{E}\left[\exp\left(p\nu_T^{1-\frac{2}{\delta^\ast}}\left(1+\int_0^T e^{\mu_t}
\bar{f}_t\mathrm{d}t\right)^{\frac{2}{\delta^\ast}}\right)\right]<+\infty
\end{equation}
with $\bar{f}_t:=f_t+\kappa_t|Y'_t|,~t\in \T$. In fact, in light of \eqref{eq:3.1} and \eqref{eq:3.8}, we have for any $p>1$,
\begin{align}\label{eq:3.11}
\begin{split}
\mathbb{E}&\left[\exp\left(p\nu_T^{1-\frac{2}{\delta^\ast}}\left(\int_0^T
e^{\mu_t}\kappa(t)|Y'_t|\mathrm{d}t\right)^{\frac{2}{\delta^\ast}}\right)\right]\\
&\leq \mathbb{E}\left[\exp\left(p\nu_T^{1-\frac{2}{\delta^\ast}}\left(\underset{t\in[0,T]}\sup (e^{\mu_t}|Y'_t|)\int_0^T\kappa(t)\mathrm{d}t\right)^{\frac{2}{\delta^\ast}}\right)\right]\\
&\leq\mathbb{E}\left[\exp\left(p \tilde K^{\frac{2}{\delta^\ast}}\nu_T^{1-\frac{2}{\delta^\ast}}\sup_{t\in\T}\left( e^{\mu_t}|Y'_t|\right)^{\frac{2}{\delta^\ast}}\right)\right]<+\infty.
\end{split}
\end{align}
The desired assertion \eqref{eq:3.10} follows immediately by \eqref{eq:3.7}, \eqref{eq:3.11} and H\"{o}lder's inequality.\vspace{0.2cm}

Next, we suppose that $g$ satisfies \eqref{eq:3.2} and $\as, g(t,Y'_t,Z'_t)\leq g'(t,Y'_t,Z'_t)$.  We will utilize the $\theta$-difference method used in \citet{BH1}. For any fixed $\theta\in(0,1)$, define
$$
\Delta_\theta U_\cdot:=\displaystyle\frac{Y_\cdot-\theta Y'_\cdot}{1-\theta}~~\mathrm{and}~~\Delta_\theta V_\cdot:=\displaystyle\frac{Z_\cdot-\theta Z'_\cdot}{1-\theta}.
$$
Then the pair $(\Delta_\theta U_\cdot,\Delta_\theta V_\cdot)$ solves the following BSDE
\begin{equation}\label{eq:3.12}
\Delta_\theta U_t=\Delta_\theta U_T+\int_{t}^{T} \Delta_\theta g(s,\Delta_\theta U_s,\Delta_\theta V_s)\mathrm{d}s-\int_{t}^{T}\Delta_\theta V_s \cdot \mathrm{d}B_s,\ \ t\in[0,T],
\end{equation}
where $\as$, for each $(y,z)\in \R\times\R^d$,
\begin{align}\label{eq:3.13}
\begin{split}
\Delta_\theta g(s,y,z):=&\frac{1}{1-\theta}[g(s,(1-\theta)y+\theta Y'_s,(1-\theta)z+\theta Y'_s)-\theta g(s,Y'_s,Z'_s) ]\\
&+\frac{\theta}{1-\theta}[g(s,Y'_s,Z'_s)-g'(s,Y'_s,Z'_s)].
\end{split}
\end{align}
It follows from \eqref{eq:3.2}, \eqref{eq:3.9} and \eqref{eq:3.13} that $\as$, for each $(y,z)\in \R\times\R^d$,
\begin{align}\label{eq:3.14}
\begin{split}
\mathbf{1}_{\{y>0\}}\Delta_\theta g(s,y,z)&\leq f_s(\omega)+\kappa_t(\omega)|Y'_s|+\beta_s(\omega)y^+ +\gamma_s(\omega)|z|^\delta\\
&=\bar{f}_s(\omega)+\beta_s(\omega) y^+ +\gamma_s(\omega)|z|^\delta.
\end{split}
\end{align}
By \eqref{eq:3.8} we know that for each $p>1$,
\begin{align*}
\begin{split}
~~~~&\mathbb{E}\left[\exp\left(p\nu_T^{1-\frac{2}{\delta^\ast}}\left(\underset{t\in[0,T]}\sup e^{\mu_t}\Delta_\theta U_t^+\right)^{\frac{2}{\delta^\ast}}\right)\right]\\&\leq\mathbb{E}\left[
\exp\left(\frac{p}{1-\theta}\nu_T^{1-\frac{2}{\delta^\ast}}\left(\underset{t\in[0,T]}\sup e^{\mu_t}(|Y_t|+|Y'_t|)\right)^{\frac{2}{\delta^\ast}}\right)\right]<+\infty,
\end{split}
\end{align*}
which along with \eqref{eq:3.10} and H\"{o}lder's inequality implies that
$$
\mathbb{E}\left[\exp\left(p\nu_T^{1-\frac{2}{\delta^\ast}}\left(\underset{t\in[0,T]}\sup (e^{\mu_t}\Delta_\theta U^+_t)+(\delta^\ast)^{\frac{\delta^\ast}{2}}+\int_0^Te^{\mu_s}\bar{f}_s\mathrm{d}s
\right)^{\frac{2}{\delta^\ast}}\right)\right]<+\infty,~~~\forall p>1.
$$
Thus, in light of \cref{rmk:2.3} along with \eqref{eq:3.14} and
\begin{equation*}\Delta_\theta U_T^+=\frac{(\xi-\theta\xi')^+}{1-\theta}=\frac{[\xi-\theta\xi+\theta(\xi-\xi')]^+}{1-\theta}\leq\xi^+,
\end{equation*}
we know that for each $t\in\T$,
\begin{equation*}
\exp\left(\nu_t^{1-\frac{2}{\delta^\ast}}(e^{\mu_t}\Delta_\theta U_t^+)^{\frac{2}{\delta^\ast}}\right)\leq \mathbb{E}\left[\exp\left(\nu_T^{1-\frac{2}{\delta^\ast}}\left(e^{\mu_T}\xi^+
+(\delta^\ast)^{\frac{\delta^\ast}{2}}+\displaystyle\int_0^Te^{\mu_s} \bar f_s\mathrm{d}s
\right)^{\frac{2}{\delta^\ast}}\right)\bigg|\mathcal{F}_t\right].
\end{equation*}
Furthermore, by virtue of inequality $e^x\geq x,\ x\geq 0$ and the last inequality along with the fact of $\nu_\cdot\geq 1$, we can deduce that for each $t\in [0,T]$,
\begin{align*}
e^{\mu_t}(Y_t-\theta Y'_t)^+\leq (1-\theta)\left(\mathbb{E}
\left[\exp\left(\nu_T^{1-\frac{2}{\delta^\ast}}
\left(e^{\mu_T}\xi^++(\delta^\ast)^{\frac{\delta^\ast}{2}}+\displaystyle\int_0^T
e^{\mu_s} \bar f_s\mathrm{d}s\right)^{\frac{2}{\delta^\ast}}\right)
\bigg|\mathcal{F}_t\right]\right)^{\frac{\delta^\ast}{2}}.
\end{align*}
Thus, the desire conclusion follows immediately by sending $\theta\To 1$ in last inequality.\vspace{0.2cm}

Moreover, if the process $\kappa_\cdot$ in \ref{A:A3} vanishes, then \eqref{eq:3.10} is just \eqref{eq:3.7}. Thus, in light of \eqref{eq:3.9*} and \eqref{eq:3.7}, we know that the conditions required in \cref{rmk:2.3} are also satisfied and then the preceding argument remains be valid. Hence, the desired conclusion holds still.\vspace{0.2cm}

The other cases can be proved in the same way as above. The proof is then complete.
\end{proof}

\def\cprime{$'$}

\end{document}